\documentclass[12pt]{article}
\usepackage{graphicx}
\usepackage{array,color}
\usepackage{amsbsy}
\usepackage{amsmath}
\usepackage{amssymb} 
\usepackage{latexsym}
\usepackage{a4wide}
\usepackage[numbers]{natbib}
\usepackage{mathrsfs}
\usepackage{bm}

\newcommand{\bqn}{\begin{eqnarray*}}
\newcommand{\eqn}{\end{eqnarray*}}
\newcommand{\um}{{\underline{m}}}
\newcommand{\cC}{{\mathcal C}}
\newcommand{\bX}{\mathbf{X}}

\newcommand{\bF}{\mathbf{F}}
\newcommand{\bSi}{\boldsymbol{\Sigma}}

\newcommand{\bbS}{{\mathbf{S}}}
\newcommand{\bbT}{{\mathbf{T}}}

\newcommand{\bbB}{{\mathbf{B}}}
\newcommand{\bbI}{{\mathbf{I}}}
\newcommand{\bbF}{{\mathbf{F}}}

\newcommand{\bxi}{{\bm\xi}}
\newcommand{\bet}{{\bm\eta}}

\newcommand{\gS}{{\bm \Sigma}}
\newcommand{\rE}{{\rm E}}

\newcommand{\ointctrclockwise}{\oint}
\newcommand{\eprf}{\hspace*{\fill}$\blacksquare$}
\numberwithin{equation}{section}  

\newtheorem{thm}{Theorem}[section]
\newtheorem{prop}{Proposition}[section]
\newtheorem{lem}{Lemma}[section]
\newtheorem{rem}{Remark}[section]

\begin{document}

\newcommand\gai[1]{{\color{blue}#1}}  

%
\begin{center}
 {\bf\Large CLT for large dimensional general Fisher matrices and its applications in high-dimensional data analysis}\\
\end{center}
\begin{center}
Shurong Zheng$^a$, Zhidong Bai$^{a}$\footnote{Corresponding email:
baizd@nenu.edu.cn } and Jianfeng Yao$^b$
\end{center}
\begin{center}
$^a$KLAS and School of Mathematics and Statistics, Northeast
Normal University, China\\
$^b$Department of Statistics and Actuarial Science, The University of Hong Kong, China
\end{center}

\centerline{\bf Abstract}

Random Fisher matrices arise naturally in multivariate statistical
analysis and  understanding the properties of its eigenvalues is
of primary importance for many  hypothesis testing problems like
testing the equality between two multivariate population
covariance matrices, or testing the independence between
sub-groups of a multivariate random vector.  This paper is
concerned with the properties of a large-dimensional Fisher matrix
when the dimension of the population is proportionally large
compared to the sample size. Most of existing works on Fisher
matrices deal with a particular Fisher matrix where populations
have i.i.d components so that the population covariance matrices
are all identity. In this paper, we consider general Fisher
matrices with arbitrary population covariance matrices. The first
main result of the paper establishes the limiting distribution of
the eigenvalues of a Fisher matrix while in a second main result,
we provide a central limit theorem for a wide class of functionals
of its eigenvalues. Some applications of these results are also
proposed for  testing hypotheses on high-dimensional covariance
matrices.


\par

\section{Introduction} \label{sec1}

  For testing the equality of variances from two populations, a
  well-known statistic is  the Fisher statistic
  defined as  the ratio of two sample
  variances.  Its multivariate counter-part is a random
  {\em Fisher matrix} defined by

\begin{equation}
  \bF:={\bf B}_1{\bf B}_2^{-1}\label{d2}
\end{equation}
where the $\bbB_j$'s are sample covariance matrices from two
independent samples, say
$\{\bxi_k, 1\le k\le n_1\}$ and $\{\bet_\ell,
1\le \ell\le n_2\}$
with population covariance matrices $\gS_1$ and $\gS_2$,
respectively.  Of primary importance are the so-called
{\em linear spectral statistics} (LSS) of
the matrix $\bF$ of form
\begin{equation}
  W_{\bf n} =\sum\limits_{i=1}^pf(\lambda_i^{\bf F}), \label{eq:LSS}
\end{equation}
where  $\lambda_i^{\bf F}$s are the eigenvalues of ${\bf F}$
with the notation ${\bf n}=(n_1,n_2) $.
Fisher matrices, especially its eigenvalues,
arise in many hypothesis testing problems in
multivariate analysis. Examples include
the test of the  equality hypothesis
$\gS_1=\gS_2$ where
the likelihood ratio (LR)   statistic is
simplified to a functional of eigenvalues of a
Fisher matrix, see \citet{BJYZ09}.
In {multivariate analysis of variance} (MANOVA), the test on the
equality of means is  reduced
to a statistic depending on a Fisher matrix which is a function of
the  ``between" sum of squares and
the ``within" sum of squares (\citet[p. 346]{Anderson03}).
In multivariate linear regression,
the likelihood ratio criterion for testing linear hypotheses about
regression coefficients is expressed as a function of the eigenvalues
of a Fisher  matrix
(\citet[p. 298]{Anderson03}). To test the independence between
sub-groups of a multivariate population, the
LR statistic is a functional of a  Fisher  matrix
defined by sub-matrices of sample covariance matrices
(\citet[p. 381]{Anderson03}).  Fisher matrices appear also in
the  canonical correlation analysis, see
\citet{YP12} for a recent account.

This paper concerns the  high-dimensional situation where the
population dimension $p$ is large compared to the sample sizes
$n_1$ and $n_2$. It is now well understood that classical procedures
as those presented in  \citet{Anderson03} become  impracticable or
dramatically lose efficiency with high-dimensional data.
For example, the deficiency of the
Hotelling's $T^2$ statistic has been reported in  \citet{D58}
and \citet{BaiSar96}.  Regarding hypothesis testing on
high-dimensional covariance matrices,
many recent works appeared in the literature, see e.g.
\cite{BJYZ09},
\cite{C10b}, \cite{Ledoit2002},
 \cite{Srivastava05}, \cite{Sriv11},
\cite{Schott07b},  and
\cite{Schott07}.
However, most of these works concern the one-sample situation and
in those treating the two-sample situation (except \cite{BJYZ09}),
the  test statistics are
often proposed though an ad-hoc distance measure so that they do not
involve the corresponding Fisher matrices.  Indeed, as it can be seen
from the multivariate analysis examples discussed earlier, Fisher
matrices and its eigenvalues appear naturally in procedures based on
the Gaussian likelihood functions.

In the literature from random matrix theory and assuming that the
dimension grows to infinity proportionally to sample  sizes,
the convergence of the eigenvalues of a  Fisher matrix
to a limiting distribution
has  been studied by several authors, see e.g.
\cite{Bai84}, \cite{BYK86}, \cite{Pill67}, \cite{PF84},
 \cite{Silv85}, and \cite{YBK83}.
As for central limit theorems for
{linear spectral statistics},
\citet{Chat09} establishes the existence of   a Gaussian
limit assuming that the populations are Gaussian.
However,  his  method doesn't provide explicit
formula for the asymptotic mean and asymptotic covariances
of the Gaussian limit.
A closely related piece of work is that of
\citet{BS04} which establishes a CLT for spectral statistics
of a general sample
covariance matrices of form ${\bf B}_1{\bf T}_p$ where ${\bf B}_1$
is a sample covariance matrix and ${\bf T}_p$ is a non-random
Hermitian matrix.  This  CLT is later refined in \cite{PanZhou08}
where the original restriction on the values of the  fourth moments of
the population components is removed.
However,  the CLT in  \cite{BS04} cannot cover spectral statistics of
a Fisher matrix  by replacing ${\bf T}_p$ by ${\bf B}_2^{-1}$
for the reason that
the centering term of this CLT would become a random
term without an explicit expression.
To overcome this difficulty, \citet{Zheng2012} establishes
a CLT for spectral statistics of a Fisher
matrix which has a non-random and explicit
centering term. In particular, the components of the
observations  ${\bxi_i}$ and ${\bet_j}$ can have arbitrary values of
the fourth moment.  To our best knowledge, this is the only CLT
reported  in the literature for spectral statistics of Fisher matrix.
However, this CLT has a severe limitation in that it is assumed that
the population covariance matrices are  equal
i.e.
$\bSi_1=\bSi_2$.
Although the derivation of this CLT is complex and highly non trivial,
it has  a small impact on the statistical problems mentioned above
where the population covariance matrices $\bSi_i$ can be arbitrary and
not necessarily equal.  Specifically for the test
of the equality hypothesis, ``$\bSi_1=\bSi_2$'' and assuming that the
population are Gaussian, this CLT enables us to find the distribution
of the LR statistic under the null hypothesis, but not under any
alternative hypothesis, that is, the size of the test can be found by
this CLT and not the power function.

The main contribution of the paper is the establishment of the central limit theorem for
linear spectral statistics  $\{ W_{\bm n} \}$
of a general Fisher matrix where
the population covariance matrices $\bSi_i$ are  arbitrary.
Under this scheme and as a preparatory step, we also establish a
limiting distribution for its eigenvalues and give an
explicit equation satisfied by its  Stieltjes transform.
Due to the fact that the population covariance matrices are arbitrary,
the establishment of these results have required several new
techniques compared to the existing literature on the central limit
theory although the general scheme remains similar to the one used in
\cite{BS04,Zheng2012}.  A significantly different tool used here is
another  CLT reported in \cite{ZBY13} for random matrix of type
${\bf S}^{-1}T$ where $\bf S$ is a standard sample covariance matrix (with
i.i.d. standardised components) and $T$ a nonnegative definite
and deterministic Hermitian matrix.
These two papers are related each other but focus on different random
matrices.

The paper is organized as follows. In Section~\ref{sec:Fmat} we
first  introduce the asymptotic scheme
 and  the technical assumptions  used, and then establishes
the limiting spectral distribution of the eigenvalues.
Section~\ref{sec:results} presents the CLT for linear spectral
statistics of general Fisher  matrices which is the main result of
the paper. Section~\ref{Computation} gives two algorithms to
approximate the limiting spectral density, mean function and
covariance function in CLT for linear spectral statistics. In
Section~\ref{example}, we discuss some applications of the results
to hypothesis testing and confidence intervals about
high-dimensional covariance matrices. Technical lemmas and  proofs
are postponed to Appendix \ref{sec:proofs}.

\section{Limiting spectral distribution of large dimensional general $F$-matrices}
\label{sec:Fmat}

Following \citet{BS04} and \citet{Zheng2012}, we will impose
the following structure on the observation model.
Assume  that the samples can be expressed as
\[ \bxi_k = \gS_1^{1/2}  {\bf X}_{\cdot k},  ~~1\le k\le n_1~;
\quad
\bet_\ell = \gS_2^{1/2}  {\bf Y}_{\cdot \ell},  ~~1\le \ell\le n_2~;
\]
where the observations matrices
\begin{eqnarray*}
  {\bf X} &:=& ({\bf X}_{\cdot 1},\cdots,{\bf X}_{\cdot n_1})
=  (X_{jk}:~ 1\le j\le p,~ 1\le k\le n_1)~,
  \\
    {\bf Y} &:=&({\bf Y}_{\cdot 1},\cdots,{\bf Y}_{\cdot n_2})=
    (Y_{j\ell}:~ 1\le j\le p,~ 1\le \ell\le n_2)~,
\end{eqnarray*}
are the upper-left corners, of size $p\times n_1$ and $p\times n_2$,
of
two independent  arrays of independent random variables
$\{X_{jk}, j,k=1,2,\cdots\}$ and $\{Y_{jk}, j,k=1,2,\cdots\}$,
respectively.   The corresponding sample covariance matrices become
\begin{eqnarray}
  {\bf B}_1&=&
 \frac{1}{n_1} \sum_{k=1}^{n_1} \bxi\bxi^* =
 \gS_1^{\frac12}\bbS_1(\gS_1^{\frac12})^*, \quad
 \text{with~~} \bbS_1 =  \frac{1}{n_1} \sum_{k=1}^{n_1}
 {\bf X}_{\cdot k}   {\bf X}_{\cdot k}  ^*  ~,\\
  {\bf B}_2&=&
 \frac{1}{n_2} \sum_{\ell=1}^{n_1} \bet\bet^* =
 \gS_2^{\frac12}\bbS_2(\gS_2^{\frac12})^*, \quad
 \text{with~~} \bbS_2 =  \frac{1}{n_2} \sum_{\ell=1}^{n_2}
 {\bf Y}_{\cdot \ell}   {\bf Y}_{\cdot \ell}  ^*  ~.
\end{eqnarray}
Because ${\bf F}= {\bf B}_1{\bf B}_2^{-1}$ has the same
eigenvalues as
 ${\bf S}_1(\bbT_p^{1/2})^*{\bf S}_2^{-1}\bbT_p^{1/2}$
where
$\bbT_p^{1/2}=\bm{\Sigma_2}^{-\frac12}\bm{\Sigma_1}^{\frac12} $,
we can define as well the Fisher matrix to be $
  \bF:={\bf S}_1(\bbT_p^{1/2})^*{\bf S}_2^{-1}\bbT_p^{1/2}
$.
It is also noticed that
obviously, the matrix ${\bf S}_2$ should be invertible
(almost surely)
so that in our asymptotic analysis, we will assume
$n_2>p$ for large $p$ and $n_2$.

Throughout the paper,
  {\em empirical spectral distribution} (or ESD) of
  square matrix refers to the empirical distribution generated by its
  eigenvalues.
We consider the following assumptions.

\begin{description}
  \item {\bf Assumption [A]} \quad The two double arrays
    $\{X_{ki}, i,k=1,2,\cdots\}$ and
    $\{Y_{ki}, i,k=1,2,\cdots\}$ consist of independent but not necessarily identically distributed random variables with mean 0 and variance 1.

  \item {\bf Assumption [B1]} \quad  For any fixed $\eta>0$ and  when $n_1, n_2,p\to\infty$,
   \begin{equation}
      \frac{1}{n_1p}\sum_{j=1}^p \sum_{k=1}^{n_1}
      \rE\left[|X_{jk}|^2
      I_{\{|X_{jk}|\geq \eta \sqrt{n_1}\}}\right]
      \rightarrow 0,
      \quad
      \frac{1}{n_2p}\sum_{j=1}^p \sum_{k=1}^{n_2}
      \rE\left[|Y_{jk}|^2
      I_{\{|Y_{jk}|\geq \eta \sqrt{n_2}\}}\right]
      \rightarrow 0.
      \label{eq:lindberg2}
    \end{equation}

  \item {\bf Assumption [B2]} \quad The two arrays are either both
    real, we then set the indicator $\kappa=2$;
    or both complex, we then set $\kappa=1$,
    with homogeneous 4th moments:
  $\rE|X_{jk}|^4=1+\kappa+\beta_x+o(1)$, $\rE|Y_{jk}|^4=1+\kappa+\beta_y+o(1)$.
    Moreover,
    for any fixed $\eta>0$ when $n_1, n_2, p\to\infty$,
   \begin{equation}
      \frac{1}{n_1p}\sum_{j=1}^p \sum_{k=1}^{n_1}
      \rE\left[|X_{jk}|^4
      I_{\{|X_{jk}|\geq \eta \sqrt{n_1}\}}\right]
      \rightarrow 0,
      \quad
      \frac{1}{n_2p}\sum_{j=1}^p \sum_{k=1}^{n_2}
      \rE\left[|Y_{jk}|^4
      I_{\{|Y_{jk}|\geq \eta \sqrt{n_2}\}}\right]
      \rightarrow 0.
      \label{eq:lindberg4}
    \end{equation}

   In addition, $\rE X_{jk}^2=o(n_1^{-1}), \rE Y_{jk}^2=o(n_2^{-1})$ when both arrays
   $\{X_{jk}\}$ and $\{Y_{jk}\}$ are complex.

  \item {\bf Assumption [C]} \quad
    The sample sizes $n_1,~n_2$ and the dimension $p$
    grow to infinity in such a way that
    \begin{equation}\label{eq:y}
      y_{n_1}:=p/n_1\rightarrow y_1\in (0, +\infty), \qquad
      \quad y_{n_2}:=p/n_2\rightarrow y_2\in (0, 1)~.
    \end{equation}
\item {\bf Assumption [D]} \quad
  The
  matrices $\bbT_p$  are  non-random and nonnegative definite
  Hermitian matrices and the sequence $\{\bbT_p\}$ is
  bounded in spectral norm. Moreover,  the
  ESD $H_p$ of $\bbT_p$ tends to a proper nonrandom
  probability measure $H$ when $p\to\infty$.
\end{description}

The assumptions \eqref{eq:lindberg2} and \eqref{eq:lindberg4} are
standard Lindeberg type conditions which are necessary  for the
existence of the limiting spectral distribution for $\bbF$, and
for the CLT for LSS of $\bbF$, respectively. Moreover, under these
conditions, the variables $X_{ik}$ and $Y_{ik}$'s  can be
truncated at size $\eta_p\sqrt{p}$ ($\eta_p\downarrow 0$) without
altering asymptotic  results.

The following notations are used throughout the paper:
\[
 {\bf n}=(n_1,n_2),\quad  {\bf y}_{{\bf n}} =(y_{n_1}, y_{n_2}), \quad
 {\bf y} =(y_1,y_2)~, \quad
 h^2=y_1+y_2-y_1y_2~.
\]
In the sequel,  the limiting results will be investigated under
the regime \eqref{eq:y}   that will be simply referred as ${\bf
n}\to\infty$. Some useful concepts are now recalled. The Stieltjes
transform of a  positive Borel measure $G$ on the real line is
defined by
\begin{equation}\label{e1}
  \displaystyle{m_G(z)\equiv\int\frac{1}{\lambda-z}dG(\lambda),\quad\quad
    z\in\mathbb{C}^{+}}=\{z:~z\in\mathbb{C}, \Im(z)>0\}.
\end{equation}
This  transform has a natural extension to the lower-half plane
by the formula
\[
  m_G(z)=\overline{m}_G (\bar z),  \quad
  \mathrm{for~~}  z\in\mathbb{C}^{-}=\{z:~z\in\mathbb{C}, \Im(z)<0\}.
\]
In addition to ${\bf F}$, we will also need several other
matrices. Table~\ref{table}  contains the notations used in the
sequel for characteristics of these matrices:  ESD, LSD and the
associated Stieltjes transforms.

\begin{table}
  \caption{Notations for distributions and Stieltjes transforms (S.T.) of
    random matrices\label{table}}
  \begin{center}
  \begin{tabular}{c@{\qquad}c@{\qquad}c}
    Matrix  & ESD / S.T.  & LSD / S.T.  \\ \hline
    $\bbF ={\bf S}_1\{\bbT_p^{1/2}\}^*{\bf S}_2^{-1}\bbT_p^{1/2}$ &   $U_{\bf n}$ /    $m_{\bf n}$
    &   $U_{\bf y}$ /    $m_{\bf y}$
    \\
     $\bX^*\{\bbT_p^{1/2}\}^*{\bf S}_2^{-1}\bbT_p^{1/2}\bX $ &   $\underline{U}_{\bf n}$ / $\underline{m}_{\bf n}$
    &   $\underline{U}_{\bf y}$ /    $\underline{m}_{\bf y}$ \\
    $\{\bbT_p^{1/2}\}^*{\bf S}_2^{-1}\bbT_p^{1/2}$  &  $G_{n_2}$ /   &   $G_{y_2}$ /
    \\ \hline
  \end{tabular}
  \end{center}
\end{table}

The matrices $\bbF$ and $\bX^*(\bbT_p^{1/2})^*{\bf S}_2^{-1}\bbT_p^{1/2}\bX $ are
companion matrices each other sharing same non null eigenvalues so
that we have
\begin{eqnarray}
  \underline{m}_{\bf n} (z)  & = & -\frac{1-y_{n_1}}{z}+y_{n_1}m_{\bf n}(z)~,\label{eq:cgn1}\\
  \underline{m}_{\bf y}(z)   & =&  -\frac{1-y_1}{z}+y_1m_{\bf y}(z)~.  \label{eq:cgn2}
\end{eqnarray}
Furthermore, when $\bSi_1=\bSi_2$, i.e., $\bbT_p=\bbI_p$, it is
well-known that the LSD $U_{\bf y}$ of $\bbF$ and its Stieltjes
transform $m_{\bf y}(z)$ can be found on Page 79 of \citet{BS09}.
As a first result of the paper, we prove the existence of $U_{\bf
y}$ and one of its characteristics for general Fisher matrix $\bbF$
where $\bbT_p$ is a Hermitian matrix.

\begin{thm}\label{thm1} Under Assumptions [A], [B1], [C]  and [D],
\begin{enumerate}
\item[(i)] The matrix $\bbS_2^{-1}\bbT_p$  has a
  non-random LSD $G_{y_2}$. Moreover, $G_{y_2}$ is characterised by the
  fact that the transform
  \[
  m_{y_2}(z)=\int_{0}^\infty\frac{t}{1-tz} dG_{y_2}(t)~,
  \]
  where
  \[
  \underline m_{y_2}(z)=-\frac{1-y_2}{z}+y_2 m_{y_2}(z)
  \]
  is the unique solution to the equation
  \begin{equation}
    z=-\frac{1}{\underline
      m_{y_2}(z)}+y_2\int\frac{dH(t)}{t+\underline m_{y_2}(z)},
    \quad z\in \mathbb{C}^+.
    \label{eqtstR1}
  \end{equation}

 \item [(ii)]
   The  Fisher  matrix
   $\bbF= \bbS_1 (\bbT_p^{1/2})^*\bbS_2^{-1}(\bbT_p^{1/2})    $
   has a non-random LSD $U_{\bf y}$. Moreover,
   $U_{\bf y}$ is
   characterised by the  fact that
   the Stieltjes
   transform $\underline m_{\bf y}(z)$ of its companion measure
   $\underline U_{\bf y}$  is the unique solution to
   the equation
   \begin{equation}
     z=\frac{ h^2 m_0(z)}{
       y_2(-1+y_2\int\frac{m_0(z)dH(t)}{t+m_0(z)})}+\frac{y_1}{y_2}m_0(z)~,  \quad z\in \mathbb{C}^+,
     \label{eqpre0}
   \end{equation}
   where $m_0(z)=\underline m_{y_2}(-\underline m_{\bf y}(z))$.
\end{enumerate}
\end{thm}
The proof of this theorem is given in Appendix \ref{sect4}.

\begin{rem}
For a given $z\in\mathbb C^+$, the equation (\ref{eqtstR1}) has a
unique solution $m_0$ such that $\Im( m_0)<0$. Then, the Stieltjes
transform $\underline m(z)$ can be computed by substituting
$z=-\underline m_{{\bf y}}$ into Equation~\eqref{eqtstR1}, i.e.
\begin{equation}
-\underline m_{{\bf y}}(z)=-\frac1{m_0(z)}+y_2\int
\frac{dH(t)}{t+m_0(z)}. \label{eqpre01}
\end{equation}
\end{rem}

In fact by \citet{Silv95}, $\underline{m}_{\bf y}(z)$ is the
unique solution to the equation
\begin{equation}\label{Jack1}
  z=\displaystyle{-\frac{1}{\um_{\bf y}(z)}+y_1\int\frac{xdG_{y_2}(x)}{1+x\um_{\bf y}(z)}}.
\end{equation}

In the sequel,  for brevity, the notations $m_{\bf y}(z)$ and
$\underline{m}_{\bf y}(z)$ will be simplified to  $m(z)$ and
$\underline{m}(z)$,  or even to   $m$ and $\um$, respectively, if no
confusion would be possible.
We will use the notations $G_{y_{n_2}}$ that are obtained by
substituting $y_{n_2}=p/n_2$ for $y_2$ in $G_{y_{2}}$.


%
\section{CLT for LSS of large dimensional general Fisher matrices}
\label{sec:results}

As explained in Introduction, we
consider {\em linear spectral statistics} (LSS) of
$\bF$
\begin{equation}
  W_{\bf n} =p\cdot U_{\bf n}(f)=\sum_{j=1}^p f(\lambda_j^{\bf F}) ~,\label{eq:LSS}
\end{equation}
where $f$ is an analytic function and $\{\lambda_j^{\bf F}\}$ are the
eigenvalues of ${\bf F}$.
More precisely, we consider a centered version
\begin{equation}\label{eq:scaled-f}
  p \left[  U_{\bf n}(f) -  U_{{\bf y}_{{\bf n}}}(f)\right] ~.
\end{equation}
where $U_{{\bf y}_{{\bf n}}}(f)=\int f(x) dU_{{\bf y}_{{\bf n}}}(x)$,
 $U_{\bf y}(x)$ is the LSD of the Fisher matrix and $U_{{\bf y}_{\bf n}}(x)$
 is obtained by substituting ${\bf y}_{{\bf n}}=(y_{n_1},y_{n_2})$ for ${\bf y}=(y_1,y_2)$
 in $U_{{\bf y}}(x)$. Due to the exact separation theorem  (see  \citet{BS99}), for
large enough $n_j$ and $p$,  the possible point mass at the origin of
$U_{\bf n}$
will coincide exactly with that of $U_{{\bf y}_{{\bf n}}}$. Therefore,
we can restrict  the integral  \eqref{eq:scaled-f} to their
continuous components on $(0,\infty)$, i.e.
\begin{eqnarray}
  p \left[  U_{\bf n} (f) -  U_{{\bf y}_{{\bf n}}}(f)   \right]
   & =&  \sum_{j=1}^p f(\lambda_j^{\bF})I_{(\lambda_j^{\bF}>0)}
  - p \int f(x) u_{{\bf y}_{{\bf n}}}(x) dx
  \label{00}
\end{eqnarray}
 where $u_{{\bf y}_{\bf n}}(x)$ is the density on $(0,\infty)$ of $U_{{\bf y}_{\bf n}}(x)$.

Regarding the central limit theory on linear spectral statistics of
random matrices, it has been well-known
(\cite{BS04,PanZhou08,Zheng2012}) that the mean and covariance
parameters of the limiting Gaussian distribution depend on the values
of the fourth moments of the initial variables. When these moments
match the Gaussian case, i.e. $\beta_x=0$ or $\beta_y=0$ in Assumption
[B2], the limiting parameters have a simpler expression. Otherwise,
they have a more involved expression that depend on other limiting
functional of sample covariance matrices. More specifically,
if $\beta_x\ne 0$, we will need the existence of the following limits
\begin{eqnarray}
  &&  \frac1p\sum\limits_{i=1}^p{\rm  E}
  \bigg[{\bf e}_i' ({\bf T}_p^{\frac12})^* {\bf S}_2^{-\frac12}
    {\bf D}_1^{-1}{\bf S}_2^{-\frac12}{\bf T}_p^{\frac12}{\bf e}_i \nonumber \\
    &&   \qquad \times~  {\bf e}_i'
    ({\bf T}_p^{\frac12})^* {\bf S}_2^{-\frac12}
    {\bf D}_1^{-1}\left(\underline{m}(z)\{{\bf
      T}_p^{\frac12}\}^*{\bf S}_2^{-1}{\bf T}_p^{\frac12}+{\bf
      I}_p\right)^{-1}{\bf S}_2^{-\frac12}{\bf T}_p^{\frac12}{\bf
      e}_i\bigg] {\longrightarrow} ~h_{m1}(z),
  \label{lx1}   \\
  &&
  \frac{1}{n_1p}\sum\limits_{j=1}^{n_1}\sum\limits_{i=1}^p {\bf    e}_i'
   ({\bf T}_p^{\frac12})^* {\bf S}_2^{-\frac12}
   [\rE_j{\bf
       D}_j^{-1}(z_1)]{\bf S}_2^{-\frac12}{\bf T}_p^{\frac12}{\bf
     e}_i       \nonumber\\
  &&\qquad \times ~
         {\bf e}_i'
         ({\bf   T}_p^{\frac12})^* {\bf S}_2^{-\frac12}
           [\rE_j{\bf D}_j^{-1}(z_2)]{\bf
    S}_2^{-\frac12}{\bf T}_p^{\frac12}{\bf
    e}_i\stackrel{i.p.}{\longrightarrow}~h_{v1}(z_1,z_2),
  \label{lx2}
\end{eqnarray}
and if $\beta_y\ne 0$, we will need the existence of the limits
\begin{eqnarray}
  &&\label{ly1}
  \frac{1}{p}\sum\limits_{i=1}^p \rE{\bf e}_i'\left(\frac1z\bbT_p-{\bf
    S}_{2}\right)^{-1}{\bf e}_i\cdot{\bf e}_i'\left(\frac1z\bbT_p-{\bf
    S}_{2}\right)^{-1}\bbT_p\left(\frac1z\bbT_p+\frac{1}{z}\underline{m}_{y_2}\left(\frac{1}{z}\right){\bf
    I}\right)^{-1}{\bf e}_i  \nonumber \\
  &&\qquad {\longrightarrow}~h_M(z)
  ~,\\
  &&\label{ly2}
  \frac{1}{n_2p}\sum\limits_{j=1}^{n_2}\sum\limits_{i=1}^p{\bf e}_i'E_j(\frac{1}{z_1}{\bf T}_p-{\bf
    S}_{2,j})^{-1}{\bf e}_i\cdot{\bf e}_i'E_j(\frac{1}{z_2}{\bf
    T}_p-{\bf S}_{2,j})^{-1}{\bf
    e}_i\stackrel{i.p.}{\longrightarrow}~h(z_1,z_2).
\end{eqnarray}
Here
\[
  {\bf  S}_{2,j}={\bf S}_2-\frac{1}{n_2}{\bf Y}_{\cdot j}{\bf   Y}^*_{\cdot j},  ~~
  {\bf D}_j(z)=({\bf S}_2^{-\frac12}{\bf T}_p^{\frac12})
  \left({\bf S}_{1}-\frac{1}{n_1}{\bf X}_{\cdot j}{\bf X}_{\cdot    j}^*\right)
  ({\bf S}_2^{-\frac12}{\bf T}_p^{\frac12})^*-z\cdot{\bf    I}_p~,
  \]
  and ${\bf e}_i$ denotes the $i$-th vector of the canonical basis of
$\mathbb{C}^p$.

The following CLT is the main result of the paper.
\begin{thm}  \label{Thm1}
  Under the   Assumptions [A], [B2], [C] and [D],
  assume that
  the limits \eqref{lx1}-\eqref{lx2} exist whenever $\beta_x\ne 0$,
  and the limits \eqref{ly1}-\eqref{ly2} exist whenever $\beta_y\ne 0$.
  Let  $f_1,\cdots,f_s$   be $s$ functions analytic in an open domain
  of the complex plane that enclosed  the support interval
  $[c_1, c_2]$ of the continuous component of the LSD
  $U_{\bf y}$.
  Then, as ${\bf n}\to\infty$,  the random vector
  \[
  \left\{  p\left[   U_{\bf n} (f_j) -  U_{{\bf y}_{{\bf n}}}(f_j)  \right]~,\quad
  1\le j\le s \right\},
  \]
  converges weakly to
  a Gaussian vector $(X_{f_1},\cdots,X_{f_s})$ with mean function
  \begin{eqnarray}
    EX_f&=&\frac{\kappa-1}{4\pi i}\ointctrclockwise_{\cC}
    f(z)~~d\log\left(\frac{h^2}{y_2}-\frac{y_1}{y_2}\cdot\frac{\left(1-y_2\int
      \frac{m_0(z)}{t+m_0(z)}dH(t)\right)^2}{1-y_2\int
      \frac{m^2_0(z)}{(t+m_0(z))^2}dH(t)}\right)\nonumber\\
    &&-\frac{\beta_x y_1}{2\pi
      i}\cdot\oint_{\cC}\frac{z^{2}\underline{m}^{3}(z)\cdot
      h_{m1}(z)}{\frac{h^2}{y_2}-\frac{y_1}{y_2}\cdot\frac{\left(1-\int
        \frac{y_2m_0(z)}{t+m_0(z)}dH(t)\right)^2}{1-\int
        \frac{y_2m_0^2(z)}{(t+m_0(z))^2}dH(t)}}dz\nonumber\\
    &&+\frac{\kappa-1}{4\pi i}\ointctrclockwise_{\cC}
    f(z)~~d\log\left(1-y_2\int
    \frac{m^2_0(z)dH(t)}{(t+m_0(z))^{2}}\right)\nonumber\\
    &&+\frac{\beta_y y_2}{2\pi
      i}\cdot\oint\um'(z)\frac{\underline{m}^3(z)m^3_0(z)
      h_M(-\frac{1}{\underline{m}(z)})} {1-y_2\int\frac{m^2_0(z)}
      {(t+m_0(z))^2}dH(t)}dz ~, \label{asym-mean0}
  \end{eqnarray}
  and covariance function
  \begin{eqnarray}
    &&{\rm Cov}{(X_{f_i}, X_{f_j})}\nonumber\\
    &=&-\frac{\beta_xy_1}{4\pi^2}\cdot
    \oint_{\cC_1}\oint_{\cC_2}\frac{\partial^2\left[z_1z_2\underline{m}(z_1)\underline{m}(z_2)h_{v1}(z_1,z_2)\right]}
         {\partial z_1\partial z_2}dz_1dz_2\nonumber\\
         &&-\frac{\kappa}{4\pi^2}\ointctrclockwise_{\cC_1}\ointctrclockwise_{\cC_2}
         \frac{f_i(z_1)f_j(z_2)}{(m_0(z_1)-m_0(z_2))^2}dm_0(z_1)dm_0(z_2)\nonumber\\
         &&-\frac{\beta_yy_2}{4\pi^2}\oint_{\cC_1}\oint_{\cC_2}
         \frac{\um'(z_1)\um'(z_2)}{\um^2(z_1)\um^2(z_2)}\frac{\partial^2\left[\underline{m}(z_1)m_0(z_1)
             \underline{m}(z_2)m_0(z_2)
             h\left(-\frac{1}{\underline{m}(z_1)},
             -\frac{1}{\underline{m}(z_2)}\right)\right]}{\partial
           (-1/\underline{m}(z_1))\partial(-1/\underline{m}(z_2))}dz_1dz_2\quad\quad\label{asym-cov0}
  \end{eqnarray}
  where  the contours $\cC$, $\cC_1$ and $\cC_2$ all enclose the
  support of $U_{\bf y}$,  and $\cC_1$ and $\cC_2$ are disjoint.
\end{thm}

Similar to CLT's developed in \cite{BS04,Zheng2012}, all the
limiting parameters depend on contour integrals using the
associated Stieltjes transforms. Some specific examples of
calculations of such integrals can be found in these references.

We next develop an important special example where the matrices
$\{{\bf T}_p\}$ are diagonal. In this case, we find explicit expressions
for the limiting functions
$h_{M}(z)$ and $h(z_1,z_2)$. This in turn simplifies the expressions
of limiting mean and covariance functions in the CLT.

\begin{prop}   \label{Thm2}
  In addition to the assumptions of Theorem~\ref{Thm1}, assume that
  the matrices  ${\bf T}_p$'s  are diagonal.
  Then, the limits \eqref{ly1} and \eqref{ly2}  exist and equal to
    \begin{eqnarray}
    h_M(z)&=&\int\frac{t}{\left(\frac{t}{z}+\frac{1}{z}\underline{m}_{y_2}\left(\frac{1}{z}\right)\right)^3}dH(t)
    ~, \label{lim-hm}\\
    h(z_1,z_2)&=&\int\frac{1}{\left(\frac{t}{z_1}+\frac{1}{z_1}\underline{m}_{y_2}\left(\frac{1}{z_1}\right)
      \right)    \left(\frac{t}{z_2}+\frac{1}{z_2}\underline{m}_{y_2}\left(\frac{1}{z_2}\right)\right)}dH(t)
    ~. \label{lim-h}
  \end{eqnarray}
 \noindent Consequently, the same conclusions as in Theorem~\ref{Thm1} hold
  where
  the last term of $\rE X_f$ in \eqref{asym-mean0}
  is simplified  to
  \begin{eqnarray}
    &&\frac{\beta_y}{4\pi i}\oint_{\cC}
    f(z)d\left(1-y_2\int\frac{m_0^2(z)dH(t)} {(t+m_0(z))^2}\right)
    \label{asym-mean}
  \end{eqnarray}
  and the last term of ${\rm Cov}(X_{f_i}, X_{f_j})$ in
  \eqref{asym-cov0} is simplified  to
  \begin{eqnarray}
    &&-\frac{\beta_y y_2}{4\pi^2}\oint_{\cC_1}\oint_{\cC_2} f_i(z_1)f_j(z_2)
    \left[\int\frac{t^2dH(t)}{(t+m_0(z_1))^2(t+m_0(z_2))^2}\right]dm_0(z_1)dm_0(z_2).
    \label{asym-cov}
  \end{eqnarray}
  where each of the contours $\cC$, $\cC_1$ and $\cC_2$ encloses the
  support of $U_{\bf y}$
  and $\cC_1$ and $\cC_2$ are disjoint.
\end{prop}

\begin{rem}
  The contours in Theorem~\ref{Thm1} and Proposition~\ref{Thm2} are taken in the
  $z$ space. In this case, the contours can be arbitrary provided
  that they enclose the support of the LSD $U_{\bf y}$.
  Since
  the integrands  are functions of $m_0$, thus the
  integrals can be taken in the $m_0$ space using the change of
  variable $z\mapsto m_0(z)$~.
\end{rem}

\begin{rem}
  When ${\bf T}_p$ is an identity matrix, (\ref{asym-mean}) and
  (\ref{asym-cov}) are the same as (3.6) and (3.7) in Zheng (2013).
  That is, Theorem 3.2 in Zheng (2012) is a special case of
  Theorem \ref{Thm2} in this paper when ${\bf T}_p={\bf I}_p$.
\end{rem}

\section{Evaluation of the asymptotic parameters  $EX_f$,
  ${\rm Cov}$ $(X_{f_i}, X_{f_j})$ and the  limiting
  density $u_{{\bf y}}(x)$}
\label{Computation}

The pratical application of Theorem~\ref{Thm1} or
Proposition~\ref{Thm2} requires to know the limiting spectral
density $u_{{\bf y}}(x)$, the asymptotic mean $EX_f$ and
covariance function ${\rm Cov}(X_{f_i}, X_{f_j})$. In particular,
the last two functions depend on some non trivial contour
integrals. In the simple  case where ${\bf T}_p=\bbI_p$ and for
simple functions like $f(x)=x^j$ (monomials) or $f(x)=\log(x)$,
analytical results can be found exactly, see  \cite{Zheng2012}.
However, this is a \textcolor{red}{very} particular case and for
general population matrices or more complex functions $f$, such
exact results are not available. In this section, we introduce
some numerical procedures to approximate these asymptotic
parameters while deliberately placing  ourselves in the context of
practical application with real data sets. In such a situation,
the sample sizes and dimension of data $(n_1, n_2, p)$ are given
and the empirical spectral distribution $H_p$ of ${\bf
T}_p=\bm{\Sigma}_1\bm{\Sigma}_2^{-1}$ is known. In this section,
we denote the eigenvalues of ${\bf T}_p$ simply by
$\{\lambda_j^0\}$ so that
$H_p(t)=\frac1p\sum\limits_{j=1}^pI_{(\lambda_j^0\leq t)}$.
However in such a concrete application situation, the LSD $H$ is
never known and we need an estimate of $H$. A very reasonable and
widely used   estimate of $H$ is indeed  just $H_p$. Here we
assume a more general estimate of the form
\begin{equation}\label{EH}
  \widehat H(t)=\sum\limits_{j=1}^pw_jI{(\lambda_j^0\leq t)}.
\end{equation}
where $\{w_j\}$ is a family of mixing weights, i.e. $w_j\ge 0$ and
$\sum w_j=1$. This form includes $H_p$ and many other interesting
estimators of $H$, e.g. a kernel estimate.

Notice that the parameters
$u_{{\bf y}}(x)$, $EX_f$ and ${\rm Cov}(X_{f_i}, X_{f_j})$ all depend
on the Stieltjes transform  $m_0(z)$.
We first approximate this transform.

\begin{lem}\label{LemCom1}
Let $z=x_z+y_z{\bf i}$ and $m_0(z) = u_0+v_0{\bf i}$  with
corresponding real and imaginary parts.
We have
\begin{equation}\label{NLE1}
  x_z=-\frac{h^2 u_0 \left(1-y_2+y_2\sum\limits_{j=1}^p
    \frac{ w_j\lambda_j^0(\lambda_j^0+u_0)}{(\lambda_j^0+u_0)^2+v_0^2 }\right)
    -h^2 v_0 \sum\limits_{j=1}^p\frac{y_2w_j\lambda_j^0v_0}
    {(\lambda_j^0+u_0 )^2+v_0^2}}
  {y_2\left(1-y_2+\sum\limits_{j=1}^p\frac{y_2w_j\lambda_j^0(\lambda_j^0+u_0 )}{(\lambda_j^0+u_0)^2
      +v_0^2 }\right)^2
    +y_2\left(\sum\limits_{j=1}^p\frac{y_2 w_j \lambda_j^0 v_0 }{(1+\lambda_j^0 u_0 )^2
      +(\lambda_j^0)^2 v_0^2 }\right)^2}+\frac{y_1 u_0 }{y_2}~,
\end{equation}
and
\begin{equation}\label{NLE2}
  y_z=-\frac{h^2 u_0 \sum\limits_{j=1}^p\frac{y_2 w_j \lambda_j^0
      v_0 }{(\lambda_j^0+u_0 )^2
      +v_0^2}
    +h^2 v_0 \left(1-y_2+\sum\limits_{j=1}^p\frac{y_2w_j\lambda_j^0(\lambda_j^0+u_0)}{(\lambda_j^0+u_0)^2
      +v_0^2}\right)}
  {y_2\left(1-y_2+\sum\limits_{j=1}^p\frac{y_2 w_j \lambda_j^0(\lambda_j^0+u_0)}{(\lambda_j^0+u_0)^2
      +v_0^2}\right)^2
    +y_2\left(\sum\limits_{j=1}^p\frac{y_2 w_j \lambda_j^0 v_0 }{(\lambda_j^0+u_0 )^2+v_0^2}\right)^2}
  +\frac{y_1 v_0 }{y_2}.
\end{equation}
\end{lem}

The  proof of Lemma~\ref{LemCom1} is given in Appendix.
Therefore given $z=x_z+y_z{\bf i}$,
$(u_0, v_0)$ are solutions of the nonlinear equations (\ref{NLE1}) and
(\ref{NLE2}). These equations  can be easily solved
using standard computing software to get numerically the values of
$(u_0, v_0)$, i.e. of $m_0(z) $.

Next, the limiting
spectral density $u_{{\bf y}}(x)$ can be approximated
as indicated  below.

\begin{rem}\label{rem1} By (\ref{eqpre01}) and (\ref{EH}) of this paper and Theorem B.10 of \citet{BS09}, we have
  \begin{equation}\label{underlinem}
    \underline
    m(z)=\frac1{m_0(z)}-y_2\int\frac{dH(t)}{t+m_0(z)}\approx\frac1{m_0(z)}-y_2\sum\limits_{j=1}^p\frac{w_j}{\lambda_j^0+m_0(z)}
  \end{equation}
  and
  \begin{equation}\label{density}
    u_{\bf y}(x)=\frac{1}{\pi y_1}\lim\limits_{\varepsilon \to o_+
      }\Im(\underline{m}(x+\varepsilon {\bf i})).
  \end{equation}
\end{rem}

\begin{rem} \label{rem2}
  The limiting functions $h_{m1}(z)$ and $h_{v1}(z_1,z_2)$ can be
  approximated as
  follows
  \begin{equation}\label{mhvh}
    \hat{h}_{m1}(z)=\frac{1}{n_1p}\sum\limits_{j=1}^{n_1}\sum\limits_{i=1}^p
    B_{1ij}(z)B_{2ij}(z),\quad \hat{h}_{v1}(z_1,
    z_2)=\frac{1}{n_1p}\sum\limits_{j=1}^{n_1}\sum\limits_{i=1}^p
    B_{1ij}(z_1)B_{1ij}(z_2)
  \end{equation}
  where
  $$
  B_{1ij}(z)={\bf e}_i'\{{\bf S}_2^{-\frac12}{\bf
    T}_p^{\frac12}\}^*{\bf D}_j^{-1}(z){\bf S}_2^{-\frac12}{\bf
    T}_p^{\frac12}{\bf e}_i
  $$
  and
  $$
  B_{2ij}(z)={\bf e}_i'\{{\bf S}_2^{-\frac12}{\bf
    T}_p^{\frac12}\}^*{\bf D}_j^{-1}(z)\left(\underline{m}(z)\{{\bf
    T}_p^{\frac12}\}^*{\bf S}_2^{-1}{\bf T}_p^{\frac12}+{\bf
    I}_p\right)^{-1}{\bf S}_2^{-\frac12}{\bf T}_p^{\frac12}{\bf e}_i
  $$
  with ${\bf D}_j(z)=({\bf S}_2^{-\frac12}{\bf T}_p^{\frac12})
  \left({\bf S}_{1}-\frac{1}{n_1}{\bf X}_{\cdot j}{\bf X}_{\cdot
    j}^*\right)({\bf S}_2^{-\frac12}{\bf T}_p^{\frac12})^*-z\cdot{\bf
    I}_p$.
\end{rem}

The following remark will give a simplified form of the asymptotic
mean function $\rE X_f$ and asymptotic covariance function ${\rm
  Cov}(X_{f_i}, X_{f_j})$.

\begin{rem} \label{LemCom4}In Theorem
    \ref{Thm2}\footnote{To be generalized to the case of Theorem 3.1.}, the mean and
  covariance functions have alternate expressions
  \begin{eqnarray}
    \rE X_f&=&-\frac{\kappa-1}{4\pi i}\ointctrclockwise_{\cC}
    f'(z)\log\left(\frac{h^2}{y_2}-\frac{y_1}{y_2}\cdot\frac{\left(1-y_2\int
      \frac{m_0(z)}{t+m_0(z)}dH(t)\right)^2}{1-y_2\int
      \frac{m_0^2(z)}{(t+m_0(z))^2}dH(t)}\right)dz\nonumber\\
    &&-\frac{\beta_x y_1}{2\pi i}\oint
    f(z)\cdot\frac{z^2\underline{m}^{3}(z) h_{m1}(z)}
    {\frac{h^2}{y_2}-\frac{y_1}{y_2}\cdot\frac{\left(1-\int
        \frac{y_2m_0(z)}{t+m_0(z)}dH(t)\right)^2}{1-\int
        \frac{y_2m_0^2(z)}{(t+m_0(z))^2}dH(t)}}dz\nonumber\\
    &&-\frac{\kappa-1}{4\pi i}\ointctrclockwise_{\cC}
    f'(z)\log\left(1-y_2\int
    \frac{m^2_0(z)dH(t)}{(t+m_0(z))^{2}}\right)dz\nonumber\\
    &&+\frac{\beta_yy_2}{4\pi i}\oint_{\cC}
    f'(z)\left(\int\frac{m_0^2(z)dH(t)} {(t+m_0(z))^2}\right)dz
    \label{ez4}
  \end{eqnarray}
  and covariance functions
  \begin{eqnarray}
    &&{\rm Cov}(X_{f_i},
    X_{f_j})\nonumber\\
    &=&-\frac{\kappa}{4\pi^2}\ointctrclockwise_{\cC_1}\ointctrclockwise_{\cC_2}
    f_i'(z_1)f_j'(z_2)\log(m_0(z_1)-m_0(z_2))dz_1dz_2\nonumber\\
    &&-\frac{\beta_y y_2}{4\pi^2}\oint\oint f_i'(z_1)f_j'(z_2)
    \left[\int\frac{t^2dH(t)}{(t+m_0(z_1))(t+m_0(z_2))}\right]dz_1dz_2\nonumber\\
    &&-\frac{\beta_x y_1}{4\pi^2}\oint\oint f_i'(z_1)f_j'(z_2)\cdot
    \left[z_1z_2\underline{m}(z_1)\underline{m}(z_2)\cdot
      h_{v1}(z_1,z_2)\right]dz_1dz_2.\label{eez3}
  \end{eqnarray}
\end{rem}

Combining the methods devised in
  Lemma \ref{LemCom1} and
  Remark \ref{rem1}-\ref{LemCom4},
  we now describe the general procedure
  to approximate the
  limiting spectral density $u_{{\bf y}}(x)$, the asymptotic mean
  and covariance functions.

\subsubsection*{\bf Algorithm 1: approximating the  limiting spectral density $u_{{\bf y}}(x)$}

Cut the support set $[c_1, c_2]$ of the LSD of Fisher matrix $\bf F$
into a mesh  set as
\[
{\cal{A}}=\left\{z_j=x_j+\varepsilon{\bf i},
x_j=c_1+\frac{(c_2-c_1)j}{m}, ~~ j=0,\ldots,m\right\}~,
\]
where $\varepsilon$ is a  small step size, e.g.
$10^{-3}$. By (\ref{NLE1}) and (\ref{NLE2}), we obtain $m_0(z_j)$
with $z_j\in{\cal{A}}$. By (\ref{underlinem}), we obtain
$\underline{m}(z_j)$ with $z_j\in{\cal{A}}$. Then by
(\ref{density}) let
\begin{equation}\label{ED}
  u_{\bf y}(x_j) \simeq \frac{1}{\pi y_1}\Im(\underline{m}(z_j))
\end{equation}
we obtain an approximation of
the density $u_{{\bf y}}(x_j)$.

\subsubsection*{\bf Algorithm 2: approximating the  asymptotic mean function
  (\ref{ez4}) and covariance function (\ref{eez3})}

\noindent {\em Step 1.}~
Chose two disjoint contours ${\cC}_1$ and ${\cC}_2$
both enclosing the support $[c_1, c_2]$ of $u_{\bf y}$
as depicted on Figure~\ref{fig1} where
$\varepsilon$ and
$\zeta$ are small numbers, e.g.
$\varepsilon=\zeta=10^{-3}$.
%

\noindent {\em Step 2.}~
Let $m_1, m_2$ be large integers, e.g.  $10^3$.
Then $\cC_1$ and $\cC_2$ are cut into a grid set as
\begin{eqnarray*}
  {\cal{A}}_1&=&\Bigg\{
  z_k=c_1-\varepsilon+\left(\zeta-\frac{2\zeta k}{m_1}\right){\bf i},~~k=0,\ldots,m_1\\
  &&
  \quad   z_{m_1+j}=c_1-\varepsilon+\frac{(c_2-c_1+2\varepsilon)j}{m_2}-\zeta{\bf i}, ~~j=0,\ldots,m_2\\
  &&
  \quad   z_{m_1+m_2+k}=c_2+\varepsilon+\left(-\zeta+\frac{2\zeta k}{m_1}\right){\bf i},~~k=0,\ldots,m_1\\
  &&
  \quad   z_{2m_1+m_2+j}=c_2+\varepsilon-\frac{(c_2-c_1+2\varepsilon)j}{m_2}+\zeta{\bf i}, ~~j=0,\ldots,m_2\Bigg\},
\end{eqnarray*}

\begin{eqnarray*}
  {\cal{A}}_2&=&\Bigg\{
  z_k=c_1-\frac{\varepsilon}{2}+\left(\frac{\zeta}{2}-\frac{\zeta k}{m_1}\right){\bf i},k=0,\ldots,m_1\\
  &&
  \quad   z_{m_1+j}=c_1-\frac{\varepsilon}{2}+\frac{(c_2-c_1+\varepsilon)j}{m_2}-\frac{\zeta}{2}{\bf i}, j=0,\ldots,m_2\\
  &&
  \quad   z_{m_1+m_2+k}=c_2+\frac{\varepsilon}{2}+\left(-\frac{\zeta}{2}+\frac{\zeta k}{m_1}\right){\bf i},k=0,\ldots,m_1\\
  &&
  \quad   z_{2m_1+m_2+j}=c_2+\frac{\varepsilon}{2}-\frac{(c_2-c_1+\varepsilon)j}{m_2}+\frac{\zeta}{2}{\bf i},
  j=0,\ldots,m_2\Bigg\}.
\end{eqnarray*}

\noindent{\em Step 3.}~  By (\ref{NLE1}) and (\ref{NLE2}), we obtain $m_0(z_j)$. By
(\ref{underlinem}), we obtain $\underline{m}(z_j)$. Then mean
function and covariance function are approximated by
\begin{eqnarray}
EX_f&\approx&
 -\frac{\kappa-1}{4\pi}\sum\limits_{j=0}^{2m_1+2m_2+3}
\Im\left[f'(z_j)\log\left(\frac{h^2}{y_2}-\frac{y_1}{y_2}\cdot\frac{\left(1-y_2\int
\frac{m_0(z_j)}{t+m_0(z_j)}dH(t)\right)^2}{1-y_2\int
\frac{m_0^2(z_j)}{(t+m_0(z_j))^2}dH(t)}\right)(z_{j+1}-z_j)\right]\nonumber\\
&&-\frac{\beta_x y_1}{2\pi}\sum\limits_{j=0}^{2m_1+2m_2+3}
\Im\left[f(z_j)\cdot\frac{z^2\underline{m}^{3}(z_j) h_{m1}(z_j)}
{\frac{h^2}{y_2}-\frac{y_1}{y_2}\cdot\frac{\left(1-\int
\frac{y_2m_0(z_j)}{t+m_0(z_j)}dH(t)\right)^2}{1-\int
\frac{y_2m_0^2(z_j)}{(t+m_0(z_j))^2}dH(t)}}(z_{j+1}-z_j)\right]\nonumber\\
&&-\frac{\kappa-1}{4\pi}\sum\limits_{j=0}^{2m_1+2m_2+3}
\Im\left[f'(z_j)\log\left(1-y_2\int
\frac{m^2_0(z_j)dH(t)}{(t+m_0(z_j))^{2}}\right)(z_{j+1}-z_j)\right]\nonumber\\
&&+\frac{\beta_yy_2}{4\pi}\sum\limits_{j=0}^{2m_1+2m_2+3}
\Im\left[f'(z_j)\left(\int\frac{m_0^2(z_j)dH(t)}
{(t+m_0(z_j))^2}\right)(z_{j+1}-z_j)\right],\quad z_j\in
{\cal{A}}_1\label{asym-meanA}
\end{eqnarray}

\begin{eqnarray}
&&{\rm Cov}(X_{f_i},
X_{f_j})\nonumber\\
&\approx&-\frac{\kappa}{4\pi^2}\sum\limits_{j,k=0}^{2m_1+2m_2+3}
\Re\left[f_i'(z_j^1)f_j'(z_k^2)\log(m_0(z_j^1)-m_0(z_k^2))(z_{j+1}^1-z_j^1)(z_{k+1}^2-z_k^2)\right]\nonumber\\
&&-\frac{\beta_y y_2}{4\pi^2}\sum\limits_{j,k=0}^{2m_1+2m_2+3}
\Re\left[f_i'(z_j^1)f_j'(z_k^2)
\left[\int\frac{t^2dH(t)}{(t+m_0(z_j^1))(t+m_0(z_k^2))}\right](z_{j+1}^1-z_j^1)(z_{k+1}^2-z_k^2)\right]\nonumber\\
&&-\frac{\beta_x y_1}{4\pi^2}\sum\limits_{j,k=0}^{2m_1+2m_2+3}
\Re\left[f_i'(z_j^1)f_j'(z_k^2)\cdot
\left[z_j^1z_k^2\underline{m}(z_j^1)\underline{m}(z_k^2)\cdot
h_{v1}(z_j^1,z_k^2)\right](z_{j+1}^1-z_j^1)(z_{k+1}^2-z_k^2)\right],\nonumber\\
&&z_j^1\in {\cal{A}}_1, z_j^2\in {\cal{A}}_2.\label{asym-covA}
\end{eqnarray}

\section{Applications to high-dimensional statistical analysis}\label{example}

In this section, we discuss two applications of the theory
  developed in the paper to two high-dimensional statistical problems.

\subsection{Power function for testing the equality of two
  high-dimensional covariance matrices}

First we consider the two-sample test of the hypothesis that
two high-dimensional covariance matrices are equal, i.e.
\begin{equation}\label{test1}
  H_0:\bm{\Sigma}_1=\bm{\Sigma}_2\quad v.s.\quad H_1:~ \gS_1\ne
  \gS_2~.
\end{equation}
By
\citet{BJYZ09}, the likelihood ratio test statistic for
(\ref{test1}) is
$$
T_N=\sum\limits_{i=1}^p\log(y_{n_1}+y_{n_2}\lambda_i)-\sum\limits_{i=1}^p\frac{y_{n_2}}{y_{n_1}
  +y_{n_2}}\log\lambda_i-\log(y_{n_1}+y_{n_2})
$$
where $\lambda_i$'s are eigenvalues of a Fisher matrix
${\bf A}{\bf B}^{-1}$ where
$${\bf
  A}=\frac{1}{n_1-1}\sum\limits_{k=1}^{n_1}\bm{\Sigma}_1^{\frac12}
({\bf X}_{\cdot k}-\bar{{\bf X}})({\bf X}_{\cdot k}-\bar{{\bf
    X}})^T\bm{\Sigma}_1^{\frac12},~~ {\bf B}=\frac{1}{n_2-1}\sum\limits_{k=1}^{n_2}
\bm{\Sigma}_2^{\frac12}({\bf Y}_{\cdot k}-\bar{{\bf Y}})({\bf Y}_{\cdot k}-\bar{{\bf
    Y}})^T\bm{\Sigma}_2^{\frac12}.$$
As mentioned in Introduction, this two-sample test has been widely
discussed in the high-dimensional context by several authors, see
e.g.
\citet{LC12} and \citet{Schott07} which used different test
statistics.
Under $H_0$ and as ${\bf n}\rightarrow \infty$,
we have
\begin{equation}
  \widetilde{T_N}=\upsilon(f)^{-\frac{1}{2}}\left[
    T_N-p \cdot F_{y_{N_1},y_{N_2}}(f)-
    m(f)\right] \stackrel{H_0}{\Rightarrow} N \left( 0, 1\right).\label{LST}
\end{equation}
where  $N_i=n_i-1$, $y_{n_i}=\frac{p}{n_i}$, $y_{N_i}=\frac{p}{N_i}$
for $i=1,2$,
and
$F_{y_{N_1},y_{N_2}}(f)$, $m(f)$ and $\upsilon(f)$ are given
in (4.5)-(4.7) of \cite{BJYZ09} with
$f(x)=\log(y_1+y_2x)-\frac{y_2}{y_1+y_2}\log x$.
The critical  region of asymptotic   size $\alpha=0.05$ is
$$
T_N>1.64\upsilon(f)^{\frac{1}{2}}+p \cdot
F_{y_{N_1},y_{N_2}}(f)+m(f).
$$
By Theorem
\ref{Thm1} in this paper, under $H_1$ we have
\begin{equation*}
  \upsilon^1(f)^{-\frac{1}{2}}\left[
    T_N-p \cdot F^1_{y_{N_1},y_{N_2}}(f)-
    m^1(f)\right] \stackrel{H_1}{\Rightarrow}N \left( 0, 1\right),
\end{equation*}
where $m^1(f)$ and $\upsilon^1(f)$ can be approximated by
(\ref{asym-meanA}) and  (\ref{asym-covA}),  and
$F^1_{y_{N_1},y_{N_2}}(f)$  by
$$
F^1_{y_{N_1},y_{N_2}}(f)=\int\limits_{c_1}^{c_2}f(x)u_{{\bf
    y}}(x)dx\approx\frac{c_2-c_1}{10^4}\sum\limits_{j=1}^{10^4}f(x_j)u_{{\bf
    y}}(x_j),\quad x_j=c_1+\frac{(c_2-c_1)j}{10^4}
$$
and $u_{{\bf y}}(x_j)$ is computed by (\ref{ED}). Since
\begin{eqnarray*}
  &&T_N\geq 1.64\upsilon(f)^{\frac{1}{2}}+p \cdot
  F_{y_{N_1},y_{N_2}}(f)+m(f)\\
  &\Leftrightarrow&\upsilon^1(f)^{-\frac{1}{2}}\left[
    T_N-p \cdot F^1_{y_{N_1},y_{N_2}}(f)-
    m^1(f)\right]\\
  &&\geq\upsilon^1(f)^{-\frac{1}{2}}\left[
    1.64\upsilon(f)^{\frac{1}{2}}+p \cdot
    F_{y_{N_1},y_{N_2}}(f)+m(f)-p \cdot
    F^1_{y_{N_1},y_{N_2}}(f)-
    m^1(f)\right]~,
\end{eqnarray*}
 the power function of the test is
$$
1-\Phi\left(\upsilon^1(f)^{-\frac{1}{2}}\left[
  1.64\upsilon(f)^{\frac{1}{2}}+p \cdot
  F_{y_{N_1},y_{N_2}}(f)+m(f)-p \cdot
  F^1_{y_{N_1},y_{N_2}}(f)-
  m^1(f)\right]\right)~,
$$
where $\Phi(\cdot)$ is the  standardized normal distribution
function.

\subsection{Confidence interval of $\theta$ in ${\bf
    T}_p(\theta)$}

As as second application, we consider  ${\bf T}_p={\bf T}_p(\theta)$, that is,
${\bf T}_p$ is determined by parameter $\theta$ which takes values
in an  interval $[a, b]$. We are interested in the confidence
interval for the  parameter $\theta$. Then using the fact
\begin{equation*}
  \upsilon^{\theta}(f)^{-\frac{1}{2}}\left[
    T_N-p \cdot F^{\theta}_{y_{N_1},y_{N_2}}(f)-
    m^{\theta}(f)\right] \stackrel{H_1}{\Rightarrow} N \left( 0, 1\right)~,
\end{equation*}
we will give a method to determine the confidence interval of
parameter $\theta$.

First cut  $[a, b]$ as
${\cal{A}}_3=\{\theta_j=a+\frac{(b-a)j}{m},~j=0,\ldots,m\}$ where
$m$ is a large integer, e.g.  $10^4$. Giving $\theta_j$,
that is, ${\bf T}_p={\bf T}_p(\theta_j)$ and using Algorithms 1-2,
we obtain $m^{\theta_j}(f)=EX_f$, $\upsilon^{\theta_j}(f)={\rm
  Cov}(X_f, X_f)$ and $F^{\theta_j}_{y_{N_1},y_{N_2}}(f)$,
$j=0,\ldots,m$. Then the confidence interval of $\theta$ is
$[\theta_L, \theta_U]$ where
$$
\theta_L=\min\left\{\theta_j:
\upsilon^{\theta_j}(f)^{-\frac{1}{2}}\left(
T_N-p \cdot F^{\theta_j}_{y_{N_1},y_{N_2}}(f)-
m^{\theta_j}(f)\right)\leq 1.64\right\}~,
$$
and
$$
\theta_U=\max\left\{\theta_j:
\upsilon^{\theta_j}(f)^{-\frac{1}{2}}\left(
T_N-p \cdot F^{\theta_j}_{y_{N_1},y_{N_2}}(f)-
m^{\theta_j}(f)\right)\leq 1.64\right\}.
$$

\section{Concluding remarks}
\label{comments}

In this paper, we have considered a general Fisher matrix ${\bf F}$ where the
(high-dimensional) population covariance matrices $\bSi_1$ and
$\bSi_2$  can be arbitrary and
not necessarily equal.  First the limiting distribution  of its
eigenvalues  has been found. Next and more importantly,  we establish
a CLT for its linear spectral  statistics. This CLT is unavoidable in
any two-sample statistical analysis with high-dimensional data.
Besides, this CLT extends and covers the CLT of \citet{Zheng2012} which
is related to standard Fisher matrices.

An important and unsolved issue on the developed theory is about the
evaluation of the limiting mean and covariance function in the CLT.
These functions have a very complex structure depending on non-trivial
contour integrals. In the special case where the matrices
$\bSi^{-1}_2\bSi_1$ are diagonal, we have proposed some simplification
though the obtained results are still complex. In
Section~\ref{Computation}, we have devised some numerical procedures
to approximate numerically these asymptotic parameters.
The advantage of these procedures is that they depend on the observed
data only. However, the accuracy of these procedure is currently
unknown.
A precise analysis of these procedures or finding other more accurate
procedures for the approximation are certainly a valuable and
challenging question in future research.


\appendix
\section{Appendix: Proofs}\label{sec:proofs}
\subsection{ Proof of Theorem \protect\ref{thm1}}\label{sect4}
Let
$$
s_{n_2}(z)=\int_0^\infty\frac{1}{t-z}dG_{n_2}(t),\quad
s_{y_2}(z)=\int_0^\infty\frac{1}{t-z}dG_{y_2}(t),
$$
be the Stieltjes transforms of the ESD and LSD $G_{y_2}(t)$ of
random matrix $(\bbT_p^{\frac12})^*{\bf
S}_2^{-1}\bbT_p^{\frac12}$, respectively. Let
\begin{equation}
m_{y_2}(z)=\int_0^\infty\frac{t}{1-tz}dG_{y_2}(t),\label{eqep3}
\end{equation}
which is the Stieltjes transform of the image measure of $G_{y_2}$
by the reciprocal transformation $\lambda \mapsto 1/\lambda$ on
$(0,\infty)$. It is easily checked that the Stieltjes transforms
are related as in
\begin{equation}\label{msrela}
  m_{y_2}(z)  =-\frac{1}{z}-\frac{1}{z^2}   s_{y_2}(1/z) ~.
\end{equation}

Similarly, consider  the  image measure and the associated Stieltjes transform
\begin{equation}\label{mnsnrela}
  m_{n_2}(z)=-\frac{1}{z}-\frac{1}{z^2}
  s_{n_2}(1/z),\quad
  m_{y_{n_2}}(z)=-\frac{1}{z}-\frac{1}{z^2}
  s_{y_{n_2}}(1/z).
\end{equation}
Let \begin{equation}\label{umy2} \underline
m_{y_2}(z)=-\frac{1-y_2}{z}+y_2m_{y_2}(z)~,
\end{equation}
 then by Theorem 2.1 of
\citet{ZBY13}, we have
\begin{equation} z=-\frac{1}{\underline
m_{y_2}(z)}+y_2\int\frac{dH(t)}{t+\underline m_{y_2}(z)},
\label{eqtst1}
\end{equation}
where $H(t)$ is the LSD of $\bbT_p$. In fact, we have
$$\underline{m}_{y_2}(z)=-\frac{1}{z}-\frac{y_2}{z^2}s_{y_2}(1/z)\quad\mbox{or}\quad
-\frac{1}{z}\underline{m}_{y_2}(\frac{1}{z})=1+y_2zs_{y_2}(z).$$
By \citet{SC95}, we have
\begin{equation}\label{um}
z=-\frac{1}{\underline
m(z)}+y_1\int\frac{tdG_{y_2}(t)}{1+t\underline
m(z)}=-\frac1{\underline m(z)}+y_1m_{y_2}(-\underline m(z)).
\end{equation}
So by (\ref{umy2}) the above equation reduces to
\begin{equation}
z=-\frac{h^2}{\underline m(z)\cdot y_2}+\frac{y_1}{y_2}\underline
m_{y_2}(-\underline m(z)). \label{eqpre2}
\end{equation}
where $h^2=y_1+y_2-y_1y_2$. Write $m_0(z)=\underline
m_{y_2}(-\underline
m(z))=\frac{1-y_2}{\underline{m}(z)}+y_2\int\frac{tdG_{y_2}(t)}{1+t\underline
m(z)}$. Replacing $z$ by $-\underline m(z)$, Eq.~(\ref{eqtst1})
becomes
\begin{equation}
-\underline m(z)=-\frac1{m_0(z)}+y_2\int \frac{dH(t)}{t+m_0(z)}.
\label{eq85}
\end{equation}
Therefore, Eq.~(\ref{eqpre2}) reduces to
\begin{equation}
z=\frac{ h^2 m_0(z)}{
y_2(-1+y_2\int\frac{m_0(z)dH(t)}{t+m_0(z)})}+\frac{y_1}{y_2}m_0(z).
\label{eqpre3}
\end{equation}
The proof of Theorem \ref{thm1} is then completed.\eprf

\subsection{Some useful identities}
\begin{lem} Let
$m_0(z)=\underline{m}_{y_2}(-\underline{m}(z))$ where
$\underline{m}(z)$ is the solution of (\ref{um}), then we have the
following identities
\begin{eqnarray}
&&1-y_1\int\frac{\underline{m}^2(z)x^2dG_{y_2}(x)}{(1+x\underline{m}(z))^2}=
\frac{h^2}{y_2}-\frac{y_1}{y_2}\cdot\frac{\left(1-\int
\frac{y_2m_0}{t+m_0}dH(t)\right)^2}{1-\int
\frac{y_2m_0^2}{(t+m_0)^2}dH(t)}\label{temp1}~,\\
&&\left[\log\left(\frac{h^2}{y_2}-\frac{y_1}{y_2}\cdot\frac{\left(1-\int
\frac{y_2m_0}{t+m_0}dH(t)\right)^2}{1-\int
\frac{y_2m_0^2}{(t+m_0)^2}dH(t)}\right)\right]'=
\frac{-2y_1\int\frac{\underline{m}^3(z)(z)x^2dG_{y_2}(x)}{(1+x\underline{m}(z))^{3}}}
{\left[1-y_1\int\frac{\underline{m}^2(z)x^2dG_{y_2}(x)}{(1+x\underline{m}(z))^{2}}\right]^2},
\label{Lem12}\\
&&\left(\frac{h^2}{y_2}-\frac{y_1}{y_2}\cdot\frac{\left(1-\int
\frac{y_2m_0}{t+m_0}dH(t)\right)^2}{1-\int
\frac{y_2m_0^2}{(t+m_0)^2}dH(t)}\right)'=
\frac{-2y_1\int\frac{\underline{m}^3(z)(z)x^2dG_{y_2}(x)}{(1+x\underline{m}(z))^{3}}}
{1-y_1\int\frac{\underline{m}^2(z)x^2dG_{y_2}(x)}{(1+x\underline{m}(z))^{2}}},
\label{Lem120}\\
&&\left[\log\left(1-y_2\int
\frac{m_0^2dH(t)}{(t+m_0)^{2}}\right)\right]'
=\frac{2\underline{m}'(z)y_2\int\frac{tm_0^3dH(t)}{(t+m_0)^{3}}}{\left(1-y_2\int
\frac{m_0^2dH(t)}{(t+m_0)^{2}}\right)^2},\label{Lem13}\\
&&m_0(z)=\frac{1}{\um(z)}\left(1-\frac{y_2}{\um(z)}s_{y_2}\left(-\frac{1}{\um(z)}\right)\right),
\quad
m_0'=\frac{-\um'm_0^2}{1-y_2\int\frac{m_0^2dH(t)}{(t+m_0)^2}}\label{Lem14}\\
&&1-y_{1}\int\frac{(\um(z))^2x^2
dG_{y_2}(x)}{(1+x\um(z))^2}=\frac{(\um(z))^2}{\um'(z)}\label{Lem15}~,\\
&&\left(1-y_2\int\frac{m_0^2(z)dH(t)}
{(t+m_0(z))^2}\right)'=2\underline{m}'(z)\frac{y_2\int\frac{m_0^3(z)t}{(t+m_0(z))^3}dH(t)}{1-y_2\int\frac{m_0^2(z)dH(t)}
{(t+m_0(z))^2}}\label{Lem16}~,
\end{eqnarray}
where $m_0'(z)=\frac{d}{dz}m_0(z)$ and
$\underline{m}'(z)=\frac{d}{dz}\underline{m}(z)$. \label{Lem1}
\end{lem}

\noindent {\bf Proof}. By (\ref{eqep3}), we have $\displaystyle
m_{y_2}'(z)=\int_0^\infty \frac{x^2dG_{y_2}(x)}{(1-xz)^2}$ where
$'$ denotes derivative.
 So by (\ref{umy2}) we have
\begin{equation}
\int\frac{x^2dG_{y_2}(x)}{(1+x\underline{m}(z))^2}=m'_{y_2}(-\underline{m}(z))=
-\frac{1-y_2}{y_2}\cdot\frac{1}{(\underline{m}(z))^2}
+\frac{1}{y_2}\cdot\underline{m}'_{y_2}(-\underline{m}(z)).\label{AE0}
\end{equation}
where
$\underline{m}_{y_2}'(-\um(z))=\frac{d}{d\xi}\underline{m}_{y_2}(\xi)_{|_{\xi=-\um(z)}}$
instead of $\frac{d}{dz}\underline{m}_{y_2}(-\um(z))$. By
(\ref{AE0}), we have
\begin{equation}
1-y_{1}\int\frac{(\um(z))^2x^2
dG_{y_2}(x)}{(1+x\um(z))^2}=\frac{h^2}{y_{2}}
-\frac{y_{1}(\um(z))^2\underline{m}_{y_2}'(-\um(z))}{y_{2}}.\label{AEt23}
\end{equation}
Differentiating both sides of (\ref{eqtst1}) and then replacing
$z$ by $-\underline m$, we obtain
\begin{equation}\label{deriva}
1=\left(\frac1{m_0^2}-y_2\int\frac{dH(t)}{(t+m_0)^2}\right)m'_{y_2}(-\underline
m).
\end{equation}
This equation, together with (\ref{eq85}), (\ref{AEt23}) and
(\ref{deriva}) imply that
\begin{equation}
1-y_1\int\frac{\underline{m}^2(z)x^2dG_{y_2}(x)}{(1+x\underline{m}(z))^2}=
\frac{h^2}{y_2}-\frac{y_1}{y_2}\cdot\frac{\left(1-\int
\frac{y_2m_0}{t+m_0}dH(t)\right)^2}{1-\int
\frac{y_2m_0^2}{(t+m_0)^2}dH(t)}. \label{eq12}
\end{equation}
Differentiating both sides of (\ref{eqpre2}) with respect to $z$,
we obtain
$$1=\frac{h^2}{y_{2}(\um(z))^2}\um'(z)-\frac{y_{1}}{y_{2}}\underline{m}_{y_2}'(-\um(z))\um'(z).$$
This implies that
$$\um'(z)=\frac{y_{2}(\um(z))^2}{h^2-y_{1}(\um(z))^2\underline{m}_{y_2}'(-\um(z))},$$
or equivalently
\begin{equation}
y_{1}(\um(z))^2\underline{m}_{y_2}'(-\um(z))=h^2-\frac{y_{2}(\um(z))^2}{\um'(z)}.
\label{AE2}
\end{equation}
 So by (\ref{AEt23}) and (\ref{AE2}), we have
\begin{equation}
  1-y_{1}\int\frac{(\um(z))^2x^2
    dG_{y_2}(x)}{(1+x\um(z))^2}=\frac{(\um(z))^2}{\um'(z)}.\label{AE23}
\end{equation}
Differentiating  both sides of (\ref{eq85}), we have
\begin{equation}\label{m0deriva}
  m_0'=\frac{-\um'm_0^2}{1-y_2\int\frac{m_0^2dH(t)}{(t+m_0)^2}}~.
\end{equation}
So we have
$$
\left(1-y_2\int\frac{m_0^2(z)dH(t)}
{(t+m_0(z))^2}\right)'=2\underline{m}'(z)\frac{y_2\int\frac{m_0^3(z)t}{(t+m_0(z))^3}dH(t)}{1-y_2\int\frac{m_0^2(z)dH(t)}
{(t+m_0(z))^2}}~.
$$
So by (\ref{eq12}) and (\ref{AE23}), we obtain
(\ref{Lem12}).
By (\ref{m0deriva}), we have the conclusion (\ref{Lem13}). By
(\ref{msrela}) and (\ref{umy2}), we have
$$
m_0(z)=\frac{1}{\um(z)}\left(1-\frac{y_2}{\um(z)}s_{y_2}\left(-\frac{1}{\um(z)}\right)\right).
$$

The proof of the lemma is completed. \eprf \vskip 0.5cm

In the sequel,  for brevity,  $s_{y_2}(z)$ will denoted as
$s(z)$ if no confusion would be possible.

\subsection{Proof of Theorem~\protect\ref{Thm1}}
\subsubsection{Deriving CLT of general Fisher matrix} Following the same techniques of truncation and
normalisation given in \citet{BS04} (see lines -9 to -6 from the
bottom of Page 559), we may assume the following additional
assumptions:
\begin{itemize}
\item $|X_{jk}|<\eta_p\sqrt{p},\ \ |Y_{jk}|<\eta_p\sqrt{p}$, \ \ for some $\eta_p\to0$, as $p\to\infty$,
\item  $EX_{jk}=0$, $EY_{jk}=0$ and $E|X_{jk}|^2=1$,
$E|Y_{jk}|^2=1$; \item $E|X_{jk}|^4=1+\kappa+\beta_x+o(1)$ and
$E|Y_{jk}|^4=1+\kappa+\beta_y+o(1)$; \item For the complex case,
$EX_{jk}^2=o(n_1^{-1})$ and $EY_{jk}^2=o(n_2^{-1})$.
\end{itemize}

We have
$$
n_1\left[\underline{m}_{\bf n}(z)-\underline{m}_{{\bf y}_{{\bf n}}}(z)\right]=
n_1\left[\underline{m}_{\bf n}(z)-\underline{m}^{\{y_{n_1},G_{n_2}\}}(z)\right]
+n_1\left[\underline{m}^{\{y_{n_1},G_{n_2}\}}(z)-\underline{m}_{{\bf y}_{{\bf n}}}(z)\right]$$
where $\underline{m}^{\{y_{n_1},G_{n_2}\}}(z)$ and
$\underline{m}_{{\bf y}_{{\bf n}}}(z)$ are the unique roots with
imaginary parts having the same signs as that of $z$ to the
following equations by (\ref{Jack1})

$$
z=-\frac{1}{\underline{m}^{\{y_{n_1},G_{n_2}\}}}+y_{n_1}\cdot\int\frac{tdG_{n_2}(t)}{1+t\underline{m}^{\{y_{n_1},G_{n_2}\}}}
\quad\mbox{and}\quad
z=-\frac{1}{\underline{m}_{{\bf y}_{{\bf n}}}}+y_{n_1}\cdot\int\frac{tdG_{y_{n_2}}(t)}
{1+t\underline{m}_{{\bf y}_{{\bf n}}}}.
$$
The proof follows two steps and  we unify the real and complex cases
with the indicator  notation $\kappa $.

\medskip \noindent{\bf Step 1.}
\quad Consider the conditional distribution of
\begin{equation}
n_1\left[\underline{m}_{\bf n}(z)-\underline{m}^{\{y_{n_1},G_{n_2}\}}(z)\right].
\label{eqstep1} \end{equation} given $\mathscr {S}_2=\big\{\mbox{all
}\bbS_{2}\big\}$. In the proof of Theorem \ref{thm1}, we have proved that $G_{n_2}$ converges to
$G_{y_2}$. Using Lemma 1.1 of Bai
and Silverstein (2004), we conclude that the conditional
distribution of
$$n_1\left[\underline{m}_{\bf n}(z)-\underline{m}^{\{y_{n_1},G_{n_2}\}}(z)\right]=
p\left[m_{\bf n}(z)-m^{\{y_{n_1},G_{n_2}\}}(z)\right]$$ given
$\mathscr S_{2}$ converges to a Gaussian process $M_1(z)$ on the
contour ${\mathcal C}$ enclosing the support $[a,b]$ of the LSD
$U_{\bf y}$ of Fisher matrix. Moreover,  its  mean function equals
\begin{eqnarray}
  \rE\left(M_1(z)|\mathscr
          {S}_2\right) &=& {(\kappa-1)\cdot\frac{
              y_1\int\underline{m}(z)^3
              x^2[1+x\underline{m}(z)]^{-3}dG_{y_{2}}(x)}{\left[1-y_1\int
                \underline{m}^2(z)x^2(1+x\underline{m}(z))^{-2}dG_{y_{2}}(x)\right]^2}}
          \nonumber \\
          && \quad
          +\beta_x\cdot\frac{h_{m1}(z)}{\left[y_1z^{2}\underline{m}^{3}(z)\right]^{-1}\cdot
            \left[1-y_1\int\frac{x^2\underline{m}^2(z)}{\{1+x\underline{m}(z)\}^2}dG_{y_2}(x)\right]}
          ~,\label{b5}
\end{eqnarray}
where we used the fact that the limit \eqref{lx1} exists for
$z\in\cC$ when $\beta_x\ne 0$ and in this case, the mean function
has then an additional term
$$
\frac{\beta_x}{p}\sum\limits_{i=1}^p\frac{\rE\left[{\bf
e}_i'\{{\bf S}_2^{-\frac12}{\bf T}_p^{\frac12}\}^*{\bf
D}_1^{-1}\{{\bf S}_2^{-\frac12}{\bf T}_p^{\frac12}\}{\bf
e}_i\cdot{\bf e}_i'\{{\bf S}_2^{-\frac12}{\bf
T}_p^{\frac12}\}^*{\bf D}_1^{-1} \left(\underline{m}(z)\{{\bf
T}_p^{\frac12}\}^*{\bf S}_2^{-1}{\bf T}_p^{\frac12}+{\bf
I}_p\right) \{{\bf S}_2^{-\frac12}{\bf T}_p^{\frac12}\}{\bf
e}_i\right]}{\left[y_1z^{2}\underline{m}^{3}(z)\right]^{-1}\cdot
\left\{1-y_1\int\frac{x^2\underline{m}^2(z)}{[1+x\underline{m}(z)]^2}dG_{y_2}(x)\right\}}.
$$
The last expression is obtained
by replacing ${\bf S}^{-1/2}$ by ${\bf S}_2^{-\frac12}{\bf T}_p^{\frac12}$ in (6.40)
of Zheng (2012).
The conditional covariance function of the process $M_1(z)$ equals
\begin{eqnarray}
  \nonumber
      {\rm Cov}(M_1(z_1), M_1(z_2)|\mathscr { S}_2)
      & = & {\kappa\cdot\left(\frac{\underline{m}'(z_1)\cdot\underline{m}'(z_2)}{(\underline{m}(z_1)-
          \underline{m}(z_2))^2}-\frac{1}{(z_1-z_2)^2}\right)}
      \nonumber \\
      \label{b6}
      && \quad +\beta_xy_1\cdot
      \frac{\partial^2\left[z_1z_2\underline{m}(z_1)\underline{m}(z_2)h_{v1}(z_1,z_2)\right]}{\partial
        z_1\partial z_2},
\end{eqnarray}
where we used the fact that the limit \eqref{lx2} exists for
$z\in\cC$ when $\beta_x\ne 0$ and in this case, the covariance
function has then an additional term obtained  by replacing ${\bf
S}^{-1/2}$ by ${\bf S}_2^{-\frac12}{\bf T}_p^{\frac12}$ in (6.41)
of Zheng (2012).

It is remarkable fact that  these limiting functions are
independent of the conditioning
$\mathscr { S}_2$, which shows that the limiting process $M_1(z)$
 is independent of the limit of the second part below.

Step 2.\quad Now, we consider the limiting process of
\begin{equation}
  n_1\left[\underline{m}^{\{y_{n_1},G_{n_2}\}}(z)
    -\underline{m}_{{\bf y}_{{\bf n}}}(z)\right]=p
  \left[m^{\{y_{n_1},G_{n_2}\}}(z)-{m}_{{{\bf y}_{{\bf n}}}}(z)\right].
  \label{eqstep2} \end{equation}
\par\noindent
By (\ref{eqep3}), we have
\begin{eqnarray}
z&=&-\frac{1}{\underline{m}_{{\bf y}_{{\bf
n}}}}+y_{n_1}\int\frac{t}{1+t\cdot\underline{m} _{{{\bf y}_{{\bf
n}}}}}dG_{y_{n_2}}(t)=-\frac{1}{\underline{m}_{{\bf y}_{{\bf n}}}}
+y_{n_1}\cdot m_{y_{n_2}}(-\underline{m}_{{\bf y}_{\bf
n}}(z)).\label{diff1}
\end{eqnarray}
On the other hand, $\underline{m}^{\{y_{n_1},G_{n_2}\}}$ is the
solution to the equation
$$
z=-\frac{1}{\underline{m}^{\{y_{n_1},G_{n_2}\}}}+y_{n_1}\int\frac{t\cdot
dG_{n_2}(t) }{1+t\cdot \underline{m}^{\{y_{n_1},G_{n_2}\}}}~,
$$
and
\begin{eqnarray}
  z&=&-\frac{1}{\underline{m}^{\{y_{n_1},G_{n_2}\}}}+y_{n_1}\int\frac{t
    dG_{n_2}(t) }{1+t\cdot \underline{m}^{\{y_{n_1},G_{n_2}\}}}\nonumber\\
  &=&-\frac{1}{\underline{m}^{\{y_{n_1},G_{n_2}\}}}+
  y_{n_1}\!\int\!\! \left\{  \frac{t dG_{n_2}(t)}{1+t    \underline{m}^{\{y_{n_1},G_{n_2}\}}}
  -  \frac{tdG_{n_2}(t) }{1+t\underline{m}_{{\bf        y}_n}}
  \right\}
  +y_{n_1}\!\int\!\!\frac{tdG_{n_2}(t) }{1+t \underline{m}_{{\bf   y}_n}} ,
  \label{diff2}
\end{eqnarray}
where
$$
\int\frac{t}{1+t\cdot\underline{m} _{{{\bf y}_{{\bf
n}}}}(z)}dG_{n_2}(t)=m_{n_2}(-\underline{m}_{{\bf y}_{{\bf
n}}}(z)).
$$
Taking the  difference of  (\ref{diff1}) and
(\ref{diff2}) yields
\begin{eqnarray*}
  0&=&-\frac{1}{\underline{m}^{\{y_{n_1},G_{n_2}\}}}+\frac{1}{\underline{m}_{{\bf y}_{{\bf n}}}}
  y_{n_1}\int \left\{  \frac{t dG_{n_2}(t)}{1+t    \underline{m}^{\{y_{n_1},G_{n_2}\}}}
  -  \frac{t dG_{n_2}(t) }{1+t\underline{m}_{{\bf        y}_n}}
  \right\}
  \\
  &&+y_{n_1}\cdot\int\frac{t\cdot dG_{n_2}(t) }{1+t\cdot
    \underline{m}_{y_n}}-y_{n_1}\cdot\int\frac{t}{1+t\cdot\underline{m}
    _{{{\bf y}_{{\bf n}}}}}dG_{y_{n_2}}(t)
\end{eqnarray*}
That is,
\begin{eqnarray*}
  0&=&\frac{\underline{m}^{\{y_{n_1},G_{n_2}\}}-\underline{m}_{{\bf
        y}_{{\bf n}}}} {\underline{m}_{{\bf y}_{{\bf
          n}}}\cdot\underline{m}^{\{y_{n_1},G_{n_2}\}}} -y_{n_1}
  \int\frac{(\underline{m}^{\{y_{n_1},G_{n_2}\}}-\underline{m}_{{\bf
        y}_{{\bf n}}})t^2dG_{n_2}(t)}
           {(1+t\underline{m}^{\{y_{n_1},G_{n_2}\}}) (1+t\underline{m}_{{\bf
                 y}_{{\bf n}}})}\\
           &&+y_{n_1}\left\}m_{n_2}(-\underline{m}_{{\bf y}_{{\bf
                 n}}})-m_{y_{n_2}}(-\underline{m}_{{\bf y}_{{\bf
                 n}}})\right\}~.
\end{eqnarray*}
Therefore, we obtain
\begin{eqnarray}
  &&n_1\cdot\left[\underline{m}^{\{y_{n_1},G_{n_2}\}}(z)-\underline{m}_{{\bf
        y}_{{\bf n}}}(z)\right]\nonumber\\
  &=&-y_{n_1}\cdot\underline{m}^{\{y_{n_1},
    G_{n_2}\}}\underline{m}_{{\bf y}_{{\bf
        n}}}\cdot\frac{n_1\left[m_{n_2}(-\underline{m}_{{\bf y}_{{\bf n}}}
      )-m_{y_{n_2}}(-\underline{m}_{{\bf y}_{{\bf n}}})\right]}
  {1-y_{n_1}\cdot\int\frac{\underline{m}_{{\bf y}_{{\bf
            n}}}\cdot\underline{m}^{\{y_{n_1},G_{n_2}\}}t^2dG_{n_2}(t)}
    {\left(1+t\underline{m}_{{\bf y}_{{\bf n}}}\right)\cdot
      \left(1+t\underline{m}^{\{y_{n_1},G_{n_2}\}}\right)}}\nonumber\\
  &=&-\underline{m}^{\{y_{n_1}, G_{n_2}\}}\underline{m}_{{\bf
      y}_{{\bf n}}}\cdot\frac{p\left[m_{n_2}(-\underline{m}_{{\bf
          y}_{{\bf n}}} )-m_{y_{n_2}}(-\underline{m}_{{\bf y}_{{\bf
            n}}})\right]} {1-y_{n_1}\cdot\int\frac{\underline{m}_{{\bf
          y}_{{\bf
            n}}}\cdot\underline{m}^{\{y_{n_1},G_{n_2}\}}t^2dG_{n_2}(t)}
    {\left(1+t\underline{m}_{{\bf y}_{{\bf n}}}\right)\cdot
      \left(1+t\underline{m}^{\{y_{n_1},G_{n_2}\}}\right)}}~. \label{tt1}
\end{eqnarray}
We then consider the limiting process of
$$p\left[m_{n_2}\left(-\underline{m}_{{\bf y}_{{\bf n}}}(z)\right)
  -m_{y_{n_2}}\left(-\underline{m}_{{\bf y}_{{\bf
        n}}}(z)\right)\right]=-\frac{p}{(\underline{m}_{{\bf y}_{{\bf
        n}}}(z))^2}\left[s_{n_2}\left(\frac{-1}{\underline{m}_{{\bf
        y}_{{\bf n}}}(z)}\right)
  -s_{y_{n_2}}\left(\frac{-1}{\underline{m}_{{\bf y}_{{\bf
          n}}}(z)}\right)\right]$$ by (\ref{mnsnrela}).
Noticing that for any $z\in {\mathbb
  C}\backslash \mathbb{R}$,
$\underline{m}_{{\bf y}_{{\bf n}}}(z)\to
\underline{m}(z)$, the limiting distribution of
$$-\frac{p}{z^2}\left[s_{n_2}\left(\frac{-1}{\underline{m}_{{\bf
        y}_{{\bf n}}}(z)}\right)
  -s_{y_{n_2}}\left(\frac{-1}{\underline{m}_{{\bf y}_{{\bf
          n}}}(z)}\right)\right]$$ is the same as that of
$$-\frac{p}{(\underline{m}_{{\bf
      y}_{{\bf
        n}}}(z))^2}\left[s_{n_2}\left(\frac{-1}{\underline{m}(z)}\right)
  -s_{y_{n_2}}\left(\frac{-1}{\underline{m}(z)}\right)\right].$$
\par\noindent
From now on, we use the notation $g(z)=-1/\underline{m}(z)$.
By Theorem 2.2 of Zheng, Bai and Yao (2013), we conclude that
$$- p g^2(z)
\left[s_{n_2}\left( g(z) \right)
-s_{y_{n_2}}\left(g(z)\right)\right]~,$$
converges weakly to a Gaussian process $M_2(\cdot)$ on $z\in
\mathcal{C}$ with mean function
\begin{eqnarray}
  \rE(M_2(z))&=&(\kappa-1)\cdot \frac{y_2\int
    \frac{t\left[1+y_2g(z)s(g(z))\right]^3dH(t)}
         {[-t\underline{m}(z)-1-y_2g(z)s(g(z))]^3}}{\left(
    1-y_2\int\frac{[1+y_2g(z)s(g(z))]^2
      dH(t)}
    {[-t\underline{m}(z)-1-y_2g(z)s(g(z))]^2}\right)^2}
  \label{b7}\\ &&
  +\frac{\beta_yy_2\left[1+y_2g(z)s(g(z))\right]^3
    h_M(g(z))}
  {1-y_2\int\frac{\left[1+y_2g(z)s(g(z))\right]^2dH(t)}
    {[-t\underline{m}(z)-1-y_2g(z)s(g(z))]^2}},\label{b7t}
\end{eqnarray}
and covariance function ${\rm Cov}(M_2(z_1), M_2(z_2))$ equaling
\begin{eqnarray}
  &&{\kappa}{g^2(z_1)g^2(z_2)}
  \left(
  \frac{\frac{\partial\left\{g(z_1)
      \left[1+y_2g(z_1)s(g(z_1))\right]\right\}}
    {\partial\left\{-1/\underline{m}(z_1)\right\}}
    \frac{\partial\left\{g(z_2)
      \left[1+y_2g(z_2)s(g(z_2))\right]\right\}}
         {\partial\left\{-1/\underline{m}(z_1)\right\}}}
       {\left\{g(z_1)\left[1+y_2g(z_1)s(g(z_1))\right]
         -g(z_2)  \left[1+y_2g(z_2)s(g(z_2))\right]\right\}^2}
       -\frac{1}{[g(z_1)- g(z_2)]^2}\right)
       \nonumber\label{b8}\\
       &&+{\beta_yy_2}{g^2(z_1)g^2(z_2)}\frac{\partial^2\left[\left(1+y_2g(z_1)
           s(g(z_1))\right)
           \left(1+y_2g(z_2)s(g(z_2))\right)h\left(g(z_1),
           g(z_2)\right)\right]}{\partial
         (-1/\underline{m}(z_1))\partial(-1/\underline{m}(z_2))}
\end{eqnarray}
for $z_1,z_2\in \mathcal{C}$, where $H(t)$ is the LSD of $\bbT_p$.
 Here we have used the fact that the limits $h_M(z)$ and
$h(z_1,z_2)$ in \eqref{ly1}-\eqref{ly2} exist whenever $\beta_y\ne
0$. Since
\[
  {1-y_{n_1}\cdot\int\frac{\underline{m}_{{\bf y}_{{\bf
            n}}}(z)\cdot\underline{m}^{\{y_{n_1},G_{n_2}\}}t^2dG_{n_2}(t)}
    {\left(1+t\underline{m}_{{\bf y}_{{\bf n}}}(z)\right)~
      \left(1+t\underline{m}^{\{y_{n_1},G_{n_2}\}}\right)}}
  \longrightarrow
      {1-y_1\int\frac{t^2\underline{m}^2(z)dG_{y_2}(t)}
    {[1+t\underline{m}(z)]^2}}~,
\]
almost surely, this limit equals $\frac{\underline{m}^2 }{ \um' }$
by (\ref{AE23}). Then by (\ref{tt1}) we have
$$
n_1\cdot\left[\underline{m}^{\{y_{n_1},G_{n_2}\}}(z)-\underline{m}_{{\bf y}_{{\bf n}}}(z)\right],
$$
converges weakly to a Gaussian process
$$
M_3(z)=-\underline{m}'(z)M_2(z),
$$
with  mean function $\rE(M_3(z))=-\underline{m}'(z)\rE M_2(z)$ and
covariance functions ${\rm Cov}(M_3(z_1),
M_3(z_2))=\underline{m}'(z_1)\underline{m}'(z_2){\rm
Cov}(M_2(z_1), M_2(z_2))$. Since the  limit process $M_1(z)$ of
$$n_1\cdot\left[\underline{m}_{\bf
    n}(z)-\underline{m}^{\{y_{n_1},G_{n_2}\}}(z)\right]$$
is independent of the ESD of $S_{n_2}$, we know that
$$\left\{n_1\cdot\left[\underline{m}_{\bf
n}(z)-\underline{m}^{\{y_{n_1},G_{n_2}\}}(z)\right],\quad
n_1\cdot\left[\underline{m}^{\{y_{n_1},G_{n_2}\}}(z)
-\underline{m}_{{\bf y}_{{\bf n}}}(z)\right]\right\}$$ converge to
a two-dimensional Gaussian process $(M_1(z),M_3(z))$ with
independent components. Consequently,
$n_1\cdot\left[\underline{m}_{\bf n}(z)-\underline{m}_{{\bf
y}_{{\bf n}}}(z)\right]$ converges weakly to  $ M_4(z)$, a
Gaussian  process with mean function
\begin{eqnarray}
  \rE(M_4(z))&=&(\kappa-1)\cdot\frac{y_1\int\underline{m}^3(z)
    x^2[1+x\underline{m}(z)]^{-3}dG_{y_2}(x)}{\left[1-y_1\int
      \underline{m}^2(z)x^2(1+x\underline{m}(z))^{-2}dG_{y_2}(x)\right]^2}\label{m1}\\
  &&+\beta_x\cdot\frac{h_{m1}(z)}{\left[y_1z^{2}\underline{m}^{3}(z)\right]^{-1}\cdot
    \left[1-y_1\int\frac{x^2\underline{m}^2(z)}{\{1+x\underline{m}(z)\}^2}dG_{y_2}(x)\right]}\label{m2}\\
  &&-(\kappa-1)\um'(z)\cdot \frac{y_2\int
    \frac{t\left[1+y_2g(z)s(g(z))\right]^3dH(t)}
         {[-t\underline{m}(z)-1-y_2g(z)s(g(z))]^3}}{\left(
    1-y_2\int\frac{[1+y_2g(z)s(g(z))]^2
      dH(t)}{[-t\underline{m}(z)-1-y_2g(z)s(g(z))]^2}\right)^2}
  \label{m3}\\
  &&-\beta_y\cdot\um'(z)\frac{y_2\left[1+y_2g(z)s(g(z))\right]^3
    h_M(g(z))}
        {1-y_2\int\frac{\left[1+y_2g(z)s(g(z))\right]^2dH(t)}
          {[-t\underline{m}(z)-1-y_2g(z)s(g(z))]^2}}~,\label{m4}
\end{eqnarray}
and covariance function
\begin{eqnarray}
  &&{\rm Cov}(M_4(z_1),M_4(z_2))\nonumber\\
  &=&\kappa\cdot\left(\frac{\underline{m}'(z_1)\cdot\underline{m}'(z_2)}{(\underline{m}(z_1)-
    \underline{m}(z_2))^2}-\frac{1}{(z_1-z_2)^2}\right)
  +\beta_xy_1\cdot\frac{\partial^2\left[z_1z_2\underline{m}(z_1)\underline{m}(z_2)h_{v1}(z_1,z_2)\right]}
  {\partial z_1\partial z_2}\nonumber\\
  &&+\kappa g'(z_1) g'(z_2)
  \frac{\frac{\partial\left\{g(z_1)
      \left[1+y_2g(z_1)s(g(z_1))\right]\right\}}
    {\partial\left\{-1/\underline{m}(z_1)\right\}}
    \frac{\partial\left\{g(z_2)
      \left[1+y_2g(z_2)s(g(z_2))\right]\right\}}
         {\partial\left\{-1/\underline{m}(z_1)\right\}}}
       {\left\{g(z_1)\left[1+y_2g(z_1)s(g(z_1))\right]
         -g(z_2)\left[1+y_2g(z_2)s(g(z_2))\right]\right\}^2}
       \nonumber\\
       &&-\kappa g'(z_1) g'(z_2)
       \frac{1}{[g(z_1)-g(z_2)]^2}\nonumber\\
       &&+  \beta_yy_2  g'(z_1) g'(z_2)
       \frac{\partial^2\left[\left(1+y_2g(z_1)
           s(g(z_1))\right)
           \left(1+y_2g(z_2)s(g(z_2))\right)h\left(g(z_1),
           g(z_2)\right)\right]}{\partial
         (-1/\underline{m}(z_1))\partial(-1/\underline{m}(z_2))}\nonumber\\
       &=&-\kappa\cdot\frac{1}{(z_1-z_2)^2}
       +\beta_xy_1\cdot\frac{\partial^2\left[z_1z_2\underline{m}(z_1)\underline{m}(z_2)h_{v1}(z_1,z_2)\right]}
       {\partial z_1\partial z_2}\quad\label{var1}\\
       &&+\kappa g'(z_1) g'(z_2)
       \frac{\frac{\partial\left\{g(z_1)
           \left[1+y_2g(z_1)s(g(z_1))\right]\right\}}
         {\partial\left\{-1/\underline{m}(z_1)\right\}}
         \frac{\partial\left\{g(z_2)
           \left[1+y_2g(z_2)s(g(z_2))\right]\right\}}
              {\partial\left\{-1/\underline{m}(z_1)\right\}}}
            {\left\{g(z_1)\left[1+y_2g(z_1)s(g(z_1))
                \right]-g(z_2)\left[1+y_2g(z_2)
                s(g(z_2))\right]\right\}^2}\quad\quad\quad\quad\label{var2}\\
            &&
            + \beta_yy_2  g'(z_1) g'(z_2)
            \frac{\partial^2\left[\left(1+y_2g(z_1)
                s(g(z_1))\right)
                \left(1+y_2g(z_2)s(g(z_2))\right)
                h\left(g(z_1),
                g(z_2)\right)\right]}{\partial
              (-1/\underline{m}(z_1))\partial(-1/\underline{m}(z_2))}. \label{var3}
\end{eqnarray}

\subsubsection{Simplifying  the mean expressions (\ref{m1}) to
  (\ref{m4}) and the covariance expressions (\ref{var2})-(\ref{var3})}
Recall that
$m_0(z)=\underline{m}_{y_2}(-\underline{m}(z))$. By (\ref{Lem12}),
we have
$$
(\ref{m1})=(\kappa-1)\cdot\frac{y_1\int\frac{\underline{m}^3(z)
    x^2}{[1+x\underline{m}(z)]^{3}}dG_{y_2}(x)}{\left[1-y_1\int
    \frac{\underline{m}^2(z)x^2}{(1+x\underline{m}(z))^{2}}dG_{y_2}(x)\right]^2}=
\frac{-(\kappa-1)}{2}\frac{d\log\left(\frac{h^2}{y_2}-\frac{y_1}{y_2}\cdot\frac{\left(1-y_2\int
    \frac{m_0(z)}{t+m_0(z)}dH(t)\right)^2}{1-y_2\int
    \frac{m_0^2(z)}{(t+m_0(z))^2}dH(t)}\right)}{dz}.
$$
By (\ref{temp1}) we have
$$
(\ref{m2})=\beta_x\cdot\frac{h_{m1}(z)}{\left[y_1z^{2}\underline{m}^{3}(z)\right]^{-1}\cdot
  \left[1-y_1\int\frac{x^2\underline{m}^2(z)}{\{1+x\underline{m}(z)\}^2}dG_{y_2}(x)\right]}
=\beta_x\cdot\frac{y_1z^{2}\underline{m}^{3}(z)\cdot
  h_{m1}(z)}{\frac{h^2}{y_2}-\frac{y_1}{y_2}\cdot\frac{\left(1-\int
    \frac{y_2m_0}{t+m_0}dH(t)\right)^2}{1-\int
    \frac{y_2m_0^2}{(t+m_0)^2}dH(t)}}~.
$$
By (\ref{Lem13}) and
(\ref{Lem14}) we have
$$
(\ref{m3})=-(\kappa-1)\um'(z)\cdot \frac{y_2\int
  \frac{t\left[1+y_2g(z)s(g(z))\right]^3dH(t)}
       {[-t\underline{m}(z)-1-y_2g(z)s(g(z))]^3}}{\left(
  1-y_2\int\frac{[1+y_2g(z)s(g(z))]^2
    dH(t)}{[-t\underline{m}(z)-1-y_2g(z)s(g(z))]^2}\right)^2}=-\frac{\kappa-1}{2}
\frac{d\log\left(1-y_2\int
  \frac{m^2_0(z)dH(t)}{(t+m_0(z))^{2}}\right)}{dz}~.
$$
We have
$$
(\ref{m4})=-\beta_y\cdot\um'(z)\frac{y_2\left[1+y_2g(z)s(g(z))\right]^3
  h_M(g(z))}
{1-y_2\int\frac{\left[1+y_2g(z)s(g(z))\right]^2dH(t)}
  {[-t\underline{m}(z)-1-y_2g(z)s(g(z))]^2}}=
-\beta_y\cdot\um'(z)\frac{y_2\underline{m}^3(z)m^3_0(z)
  h_M(g(z))} {1-y_2\int\frac{m^2_0(z)}
  {(t+m_0(z))^2}dH(t)}~.
$$
By (\ref{Lem14}) we have
\begin{eqnarray*}
  (\ref{var2})&=&\kappa g'(z_1) g'(z_2)
  \frac{\frac{\partial\left\{g(z_1)
      \left[1+y_2g(z_1)s(g(z_1))\right]\right\}}
    {\partial\left\{-1/\underline{m}(z_1)\right\}}
    \frac{\partial\left\{g(z_2)
      \left[1+y_2g(z_2)s(g(z_2))\right]\right\}}
         {\partial\left\{-1/\underline{m}(z_1)\right\}}}
       {\left\{g(z_1)\left[1+y_2g(z_1)s(g(z_1))
           \right]-g(z_2)\left[1+y_2g(z_2)
           s(g(z_2))\right]\right\}^2}\\
       &=&\kappa\cdot\frac{1}{(m_0(z_1)-m_0(z_2))^2}\frac{\partial
         m_0(z_1)}{\partial z_1}\frac{\partial m_0(z_2)}{\partial z_2},
\end{eqnarray*}
and
$$
(\ref{var3})=\beta_yy_2 g'(z_1) g'(z_2)
\frac{\partial^2\left[\underline{m}(z_1)m_0(z_1)
    \underline{m}(z_2)m_0(z_2)
    h\left(g(z_1),
    g(z_2)\right)\right]}{\partial
  (-1/\underline{m}(z_1))\partial(-1/\underline{m}(z_2))}.
$$
So we obtain
\begin{eqnarray}
  \lefteqn{-\frac{1}{2\pi i}\ointctrclockwise_{\cC}
  f_i(z)\cdot(\ref{m1})dz} \nonumber
  \\ &&=\frac{\kappa-1}{4\pi
    i}\ointctrclockwise_{\cC}
  f_i(z)~~d\log\left(\frac{h^2}{y_2}-\frac{y_1}{y_2}\cdot\frac{\left(1-y_2\int
    \frac{m_0}{t+m_0}dH(t)\right)^2}{1-y_2\int
    \frac{m_0^2}{(t+m_0)^2}dH(t)}\right)\label{m10},\\
  &&-\frac{1}{2\pi i}\ointctrclockwise_{\cC}
  f_i(z)\cdot(\ref{m2})dz=-\frac{\beta_x}{2\pi
    i}\cdot\oint\frac{y_1z^{2}\underline{m}^{3}(z)\cdot
    h_{m1}(z)}{\frac{h^2}{y_2}-\frac{y_1}{y_2}\cdot\frac{\left(1-\int
      \frac{y_2m_0}{t+m_0}dH(t)\right)^2}{1-\int
      \frac{y_2m_0^2}{(t+m_0)^2}dH(t)}}dz\label{m20},\\
  &&-\frac{1}{2\pi i}\ointctrclockwise_{\cC}
  f_i(z)\cdot(\ref{m3})dz= \frac{\kappa-1}{4\pi
    i}\ointctrclockwise_{\cC} f_i(z)~~d\log\left(1-y_2\int
  \frac{m^2_0(z)dH(t)}{(t+m_0(z))^{2}}\right)\label{m30},\\
  &&-\frac{1}{2\pi i}\ointctrclockwise_{\cC}
  f_i(z)\cdot(\ref{m4})dz=\frac{\beta_y}{2\pi
    i}\cdot\oint\um'(z)\frac{y_2\underline{m}^3(z)m^3_0(z)
    h_M(g(z))} {1-y_2\int\frac{m^2_0(z)}
    {(t+m_0(z))^2}dH(t)}dz ~,     \label{m40}
\end{eqnarray}
for the mean function, where $h^2=y_1+y_2-y_1y_2$, and
\begin{eqnarray}
  &&-\frac{1}{4\pi^2}\ointctrclockwise_{\cC_1}\ointctrclockwise_{\cC_2}
  f_i(z_1)f_j(z_2)\cdot(\ref{var1})dz\nonumber \\
  && \qquad =-\frac{\beta_xy_1}{4\pi^2}\cdot
  \oint\oint\frac{\partial^2\left[z_1z_2\underline{m}(z_1)\underline{m}(z_2)h_{v1}(z_1,z_2)\right]}
     {\partial z_1\partial z_2}dz_1dz_2~,\label{var10}\\
     &&-\frac{1}{4\pi^2}\ointctrclockwise_{\cC_1}\ointctrclockwise_{\cC_2}
     f_i(z_1)f_j(z_2)\cdot(\ref{var2})dz\nonumber\\
     &&\qquad =
     -\frac{\kappa}{4\pi^2}\ointctrclockwise_{\cC_1}\ointctrclockwise_{\cC_2}
                 \frac{f_i(z_1)f_j(z_2)}{(m_0(z_1)-m_0(z_2))^2}dm_0(z_1)dm_0(z_2)~,\label{var20}
                 \\
                 &&-\frac{1}{4\pi^2}\ointctrclockwise_{\cC_1}\ointctrclockwise_{\cC_2}
                 f_i(z_1)f_j(z_2)\cdot(\ref{var3})dz
                 =-\frac{\beta_y}{4\pi^2}\oint\oint
                 \frac{y_2\um'(z_1)\um'(z_2)}{\um^2(z_1)\um^2(z_2)}
          \nonumber \\
          &&
          \qquad \times
          \frac{\partial^2\left[\underline{m}(z_1)m_0(z_1)
              \underline{m}(z_2)m_0(z_2)
      h\left(g(z_1),
      g(z_2)\right)\right]}{\partial
    (-1/\underline{m}(z_1))\partial(-1/\underline{m}(z_2))}dz_1dz_2~, \label{var30}
\end{eqnarray}
for the covariance function.
The respective sums lead to the mean and covariance functions of the
Theorem.

\subsection{Proof of Proposition~\protect\ref{Thm2}}


When the matrices ${\bf T}_p$ are disgonal,
we first find the limit functions $h_M(z)$ and $h(z_1,z_2)$.
This will lead to the simplication of the terms  (\ref{m40})
and  (\ref{var30}).
We have
\begin{eqnarray*}
(\frac1z\bbT_p-\bbS_{2,k})^{-1}&=&(\frac1z\bbT_p-\frac{n-1}{n}E\beta_{12}(z){\bf
I})^{-1}+B(z)+E\beta_{12}(z)\cdot A(z)+C(z),
\end{eqnarray*}
where
$$
{\bf
  A}(z)=\sum\limits_{i\not=k}(\frac1z\bbT_p-\frac{n-1}{n}E\beta_{12}(z){\bf
  I})^{-1}(\bm{\alpha}_{i}\bm{\alpha}_{i}'-\frac1n)(\frac1z\bbT_p-\bbS_{ik})^{-1},
$$
$$
{\bf
  B}(z)=\sum\limits_{i\not=k}(\beta_{ik}(z)-E\beta_{12}(z))(\frac1z\bbT_p-\frac{n-1}{n}E\beta_{12}(z){\bf
  I})^{-1}\bm{\alpha}_i\bm{\alpha}_i'(\frac1z\bbT_p-\bbS_{ik})^{-1}
,$$
$$
{\bf C}(z)=\frac1n\cdot
E\beta_{12}(z)(\frac1z\bbT_p-\frac{n-1}{n}E\beta_{12}(z){\bf
  I})^{-1}\sum\limits_{i\not=k}\left[(\frac1z\bbT_p-\bbS_{ik})^{-1}-(\frac1z\bbT_p-\bbS_k)^{-1}\right],
$$
with
$\beta_{ik}(z)=\frac{1}{1-\bm{\alpha}_i'(\frac1z\bbT_p-\bbS_{2,ik})^{-1}\bm{\alpha}_i}$,
$\bbS_{2,ik}=\bbS_{2}-\bm{\alpha}_i\bm{\alpha}_i^*-\bm{\alpha}_k\bm{\alpha}_k^*$
and $\bm{\alpha}_i=\frac{1}{\sqrt{n_2}}{\bf Y}_i.$ Then we have
\begin{eqnarray*}
  {\bf e}_l'{\bf A}(z){\bf e}_l&=&\sum\limits_{i\not=k}{\bf
    e}_l'(\frac1z\bbT_p-\frac{n-1}{n}E\beta_{12}(z){\bf
    I})^{-1}(\bm{\alpha}_{i}\bm{\alpha}_{i}'-\frac1n)(\frac1z\bbT_p-\bbS_{2,ik})^{-1}{\bf
    e}_l\\
  &=&\sum\limits_{i\not=k}\bm{\alpha}_{i}'(\frac1z\bbT_p-\bbS_{2,ik})^{-1}{\bf
    e}_l{\bf e}_l'(\frac1z\bbT_p-\frac{n-1}{n}E\beta_{12}(z){\bf
    I})^{-1}\bm{\alpha}_i\\
  &&-\sum\limits_{i\not=k}\frac1n{\rm
    tr}(\frac1z\bbT_p-\bbS_{2,ik})^{-1}{\bf e}_l{\bf
    e}_l'(\frac1z\bbT_p-\frac{n-1}{n}E\beta_{12}(z){\bf
    I})^{-1}\\
  &=&\sum\limits_{i\not=k}\hat{\gamma}^1_i~,
\end{eqnarray*}
where
\begin{eqnarray*}
  \hat{\gamma}^1_i&=&\bm{\alpha}_{i}'(\frac1z\bbT_p-\bbS_{2,ik})^{-1}{\bf
    e}_l{\bf e}_l'(\frac1z\bbT_p-\frac{n-1}{n}E\beta_{12}(z){\bf
    I})^{-1}\bm{\alpha}_i\\
  &&-\sum\limits_{i\not=k}\frac1n{\rm
    tr}(\frac1z\bbT_p-\bbS_{2,ik})^{-1}{\bf e}_l{\bf
    e}_l'(\frac1z\bbT_p-\frac{n-1}{n}E\beta_{12}(z){\bf I})^{-1}~.
\end{eqnarray*}
We also have
\begin{eqnarray*}
  &&{\bf e}_l'{\bf A}(z)\bbT_p\left(\frac1z\bbT_p-(1+y_2zs(z)){\bf
    I}\right)^{-1}{\bf e}_l\\
  &=&\sum\limits_{i\not=k}{\bf
    e}_l'(\frac1z\bbT_p-\frac{n-1}{n}E\beta_{12}(z){\bf
    I})^{-1}(\bm{\alpha}_{i}\bm{\alpha}_{i}'-\frac1n)(\frac1z\bbT_p-\bbS_{2,ik})^{-1}
  \bbT_p\left(\frac1z\bbT_p-(1+y_2zs(z)){\bf I}\right)^{-1}{\bf
    e}_l\\
  &=&\sum\limits_{i\not=k}\bm{\alpha}_{i}'(\frac1z\bbT_p-\bbS_{2,ik})^{-1}
  \bbT_p\left(\frac1z\bbT_p-(1+y_2zs(z)){\bf I}\right)^{-1}{\bf
    e}_l{\bf e}_l'(\frac1z\bbT_p-\frac{n-1}{n}E\beta_{12}(z){\bf
    I})^{-1}\bm{\alpha}_{i}\\
  &&-\sum\limits_{i\not=k}\frac{1}{n}{\rm
    tr}(\frac1z\bbT_p-\bbS_{2,ik})^{-1}
  \bbT_p\left(\frac1z\bbT_p-(1+y_2zs(z)){\bf I}\right)^{-1}{\bf
    e}_l{\bf e}_l'(\frac1z\bbT_p-\frac{n-1}{n}E\beta_{12}(z){\bf
    I})^{-1}\\
  &=&\sum\limits_{i\not=k}\hat{\gamma}^2_i~,
\end{eqnarray*}
where
\begin{eqnarray*}
  \hat{\gamma}^2_i&=&\bm{\alpha}_{i}'(\frac1z\bbT_p-\bbS_{2,ik})^{-1}
  \bbT_p\left(\frac1z\bbT_p-(1+y_2zs(z)){\bf I}\right)^{-1}{\bf
    e}_l{\bf e}_l'(\frac1z\bbT_p-\frac{n-1}{n}E\beta_{12}(z){\bf
    I})^{-1}\bm{\alpha}_{i}\\
  &&-\sum\limits_{i\not=k}\frac{1}{n}{\rm
    tr}(\frac1z\bbT_p-\bbS_{2,ik})^{-1}
  \bbT_p\left(\frac1z\bbT_p-(1+y_2zs(z)){\bf I}\right)^{-1}{\bf
    e}_l{\bf e}_l'(\frac1z\bbT_p-\frac{n-1}{n}E\beta_{12}(z){\bf
    I})^{-1}.
\end{eqnarray*}
So we obtain
\begin{eqnarray*}
  &&\left|\rE{\bf e}_l'{\bf A}(z){\bf e}_l\cdot{\bf e}_l'{\bf
    A}(z)\bbT_p\left(\frac1z\bbT_p-(1+y_2zs(z)){\bf I}\right)^{-1}{\bf
    e}_l\right|^2\\
  &=&\left|\rE\sum\limits_{i\not=k}\hat{\gamma}^1_i\hat{\gamma}^2_i\right|^2=
  \sum\limits_{i\not=k}\rE\left|\hat{\gamma}^1_i\hat{\gamma}^2_i\right|^2\leq
  \sum\limits_{i\not=k}\sqrt{\rE|\hat{\gamma}^1_i|^4\rE|\hat{\gamma}^2_i|^4}\leq
  K\cdot\eta_n^4.
\end{eqnarray*}
Similarly,
$$
\left|\rE{\bf e}_l'{\bf A}(z){\bf e}_l\cdot{\bf e}_l'{\bf
  B}(z)\bbT_p\left(\frac1z\bbT_p-(1+y_2zs(z)){\bf I}\right)^{-1}{\bf
  e}_l\right|^2\leq K\cdot\eta_n^4~,
$$
$$
\left|\rE{\bf e}_l'{\bf A}(z){\bf e}_l\cdot{\bf e}_l'{\bf
  C}(z)\bbT_p\left(\frac1z\bbT_p-(1+y_2zs(z)){\bf I}\right)^{-1}{\bf
  e}_l\right|^2\leq K\cdot\eta_n^4~,
$$
$$
\left|\rE{\bf e}_l'{\bf B}(z){\bf e}_l\cdot{\bf e}_l'{\bf
  A}(z)\bbT_p\left(\frac1z\bbT_p-(1+y_2zs(z)){\bf I}\right)^{-1}{\bf
  e}_l\right|^2\leq K\cdot\eta_n^4~,
$$
$$
\left|\rE{\bf e}_l'{\bf B}(z){\bf e}_l\cdot{\bf e}_l'{\bf
  B}(z)\bbT_p\left(\frac1z\bbT_p-(1+y_2zs(z)){\bf I}\right)^{-1}{\bf
  e}_l\right|^2\leq K\cdot\eta_n^4~,
$$
$$
\left|\rE{\bf e}_l'{\bf B}(z){\bf e}_l\cdot{\bf e}_l'{\bf
  C}(z)\bbT_p\left(\frac1z\bbT_p-(1+y_2zs(z)){\bf I}\right)^{-1}{\bf
  e}_l\right|^2\leq K\cdot\eta_n^4~,
$$
$$
\left|\rE{\bf e}_l'{\bf C}(z){\bf e}_l\cdot{\bf e}_l'{\bf
  A}(z)\bbT_p\left(\frac1z\bbT_p-(1+y_2zs(z)){\bf I}\right)^{-1}{\bf
  e}_l\right|^2\leq K\cdot\eta_n^4~,
$$
$$
\left|\rE{\bf e}_l'{\bf C}(z){\bf e}_l\cdot{\bf e}_l'{\bf
  B}(z)\bbT_p\left(\frac1z\bbT_p-(1+y_2zs(z)){\bf I}\right)^{-1}{\bf
  e}_l\right|^2\leq K\cdot\eta_n^4~,
$$
$$
\left|\rE{\bf e}_l'{\bf C}(z){\bf e}_l\cdot{\bf e}_l'{\bf
  C}(z)\bbT_p\left(\frac1z\bbT_p-(1+y_2zs(z)){\bf I}\right)^{-1}{\bf
  e}_l\right|^2\leq K\cdot\eta_n^4~.
$$
Then it is easy to obtain
$$
\begin{array}{lll}
  &&\frac{1}{p}\sum\limits_{j=1}^p\rE{\bf
    e}_j'\left(\frac1z\bbT_p-{\bf S}_{2,k}\right)^{-1}{\bf e}_j\cdot
  {\bf e}_j'\left(\frac1z\bbT_p-{\bf
    S}_{2,k}\right)^{-1}\bbT_p\left(\frac1z\bbT_p-(1+y_2zs(z)){\bf
    I}\right)^{-1}{\bf e}_j\\
  &=&\frac{1}{p}\sum\limits_{j=1}^p{\bf
    e}_j'(\frac1z\bbT_p-\frac{n-1}{n}E\beta_{12}(z){\bf I})^{-1}{\bf
    e}_{j}\cdot{\bf
    e}_j'\left(\frac1z\bbT_p-\frac{n-1}{n}E\beta_{12}(z){\bf
    I}\right)^{-1}\bbT_p\left(\frac1z\bbT_p-(1+y_2zs(z)){\bf
    I}\right)^{-1}{\bf e}_j\\
  &&+o(1)\\
  &=&\frac{1}{p}\sum\limits_{j=1}^p{\bf
    e}_j'(\frac1z\bbT_p-(1+y_2zs(z)){\bf I})^{-1}{\bf e}_{j}\cdot{\bf
    e}_j'\left(\frac1z\bbT_p-(1+y_2zs(z)){\bf
    I}\right)^{-1}\bbT_p\left(\frac1z\bbT_p-(1+y_2zs(z)){\bf
    I}\right)^{-1}{\bf e}_j\\
  &&+o(1)~,
\end{array}
$$
and
\begin{eqnarray*}
  &&\frac{1}{pn_2}\sum\limits_{j=1}^p\sum\limits_{i=1}^{n_2}{\bf
    e}_j'\rE_i(\frac{1}{z_1}{\bf T}_p-{\bf S}_{2,i})^{-1}{\bf
    e}_j\cdot{\bf e}_j'\rE_i(\frac{1}{z_2}{\bf T}_p-{\bf
    S}_{2,i})^{-1}{\bf e}_j\\
  &=&\frac{1}{p}\sum\limits_{j=1}^p{\bf
    e}_j'(\frac{1}{z_1}\bbT_p-\frac{n-1}{n}E\beta_{12}(z_1){\bf
    I})^{-1}{\bf e}_j\cdot {\bf
    e}_j'(\frac{1}{z_2}\bbT_p-\frac{n-1}{n}E\beta_{12}(z_2){\bf
    I})^{-1}{\bf e}_j+o_p(1)\\
  &=&\frac{1}{p}\sum\limits_{j=1}^p{\bf
    e}_j'(\frac{1}{z_1}\bbT_p-(1+y_2z_1s(z_1)){\bf I})^{-1}{\bf
    e}_{j}\cdot{\bf e}_j'(\frac{1}{z_2}\bbT_p-(1+y_2z_2s(z_2)){\bf
    I})^{-1}{\bf e}_{j}+o_p(1).
\end{eqnarray*}
If $\bbT_p$ is diagonal, then
\begin{eqnarray*}
  {\bf e}_j'(\frac{1}{z_1}\bbT_p-(1+y_2z_1s(z_1)){\bf I})^{-1}{\bf
    e}_{j}&=&\frac{1}{\frac{\lambda_j^0}{z_1}-(1+y_2z_1s(z_1))} ~,
\end{eqnarray*}
and
\begin{eqnarray*}
  {\bf e}_j'\left(\frac1z\bbT_p-(1+y_2zs(z)){\bf
    I}\right)^{-1}\bbT_p\left(\frac1z\bbT_p-(1+y_2zs(z)){\bf
    I}\right)^{-1}{\bf e}_j
  &=&\frac{\lambda_j^0}{\left(\frac{1}{z}\lambda_j^0-(1+y_2zs(z))\right)^2},
\end{eqnarray*}
where $\lambda_j^0$s are eigenvalues of ${\bf T}_p$. So we obtain
\begin{eqnarray*}
  &&\frac{1}{p}\sum\limits_{l=1}^p\left[\left(\frac1z\bbT_p-{\bf
      S}_{2k}\right)^{-1}\right]_{ll} \left[\left(\frac1z\bbT_p-{\bf
      S}_{2k}\right)^{-1}\bbT_p\left(\frac1z\bbT_p-(1+y_2zs(z)){\bf
      I}\right)^{-1}\right]_{ll}\\
  &=&\frac{1}{p}\sum\limits_{j=1}^p\frac{\lambda_j^0}{\left(\frac{1}{z}\lambda_j^0-(1+y_2zs(z))\right)^3}+o(1)\\
  &=&\int\frac{t}{\left(\frac{t}{z}-(1+y_2zs(z))\right)^3}dH_p(t)+o(1)\\
  &=&\int\frac{t}{\left(\frac{t}{z}-(1+y_2zs(z))\right)^3}dH(t)+o(1)
\end{eqnarray*}
and
\begin{eqnarray*}
  &&\frac{1}{p}\sum\limits_{j=1}^p[(\frac{1}{z_1}{\bf T}_p-{\bf
      S}_{2i})^{-1}]_{jj} [(\frac{1}{z_2}{\bf T}_p-{\bf
      S}_{2i})^{-1}]_{jj}\\
  &=&\frac{1}{p}\sum\limits_{j=1}^p\frac{1}{\frac{\lambda_j^0}{z_1}-(1+y_2z_1s(z_1))}
  \frac{1}{\frac{\lambda_j^0}{z_2}-(1+y_2z_2s(z_2))}+o(1)\\
  &=&\int\frac{1}{\left(\frac{t}{z_1}-(1+y_2z_1s(z_1))\right)\left(\frac{t}{z_2}-(1+y_2z_2s(z_2))\right)}dH_p(t)+o(1)\\
  &=&\int\frac{1}{\left(\frac{t}{z_1}-(1+y_2z_1s(z_1))\right)\left(\frac{t}{z_2}-(1+y_2z_2s(z_2))\right)}dH(t)+o(1).
\end{eqnarray*}
That is, we have found the limits \eqref{ly1}-\eqref{ly2} with
$$
h_M(z)=\int\frac{t}{\left(\frac{t}{z}-(1+y_2zs(z))\right)^3}dH(t)~,
$$
and
$$
h(z_1,z_2)=\int\frac{1}{\left(\frac{t}{z_1}-(1+y_2z_1s(z_1))\right)\left(\frac{t}{z_2}-(1+y_2z_2s(z_2))\right)}dH(t).
$$
So we have
\begin{eqnarray*}
  h_M\left(g(z)\right)
  &=&\int\frac{t}{\left(-t\underline{m}(z)-\left(1+y_2g(z)
    s\left(g(z)\right)\right)\right)^3}dH(t)\\
  &=&-\frac{1}{\underline{m}^3(z)}\int\frac{t}{\left(t+\frac{1}{\underline{m}(z)}\left(1+y_2g(z)
    s\left(g(z)\right)\right)\right)^3}dH(t)\\
  &=&-\frac{1}{\underline{m}^3(z)}\int\frac{t}{\left(t+m_0(z)\right)^3}dH(t)~~(\mbox{by}~(\ref{Lem14})),\\
\end{eqnarray*}
and
$$
h\left(g(z_1),g(z_2)\right)=
\frac{1}{\underline{m}(z_1)\underline{m}(z_2)}\int\frac{1}{\left(t+m_0(z_1)\right)\left(t+m_0(z_2)\right)}dH(t).
$$
Then we obtain
\begin{eqnarray*}
  (\ref{m4})&=&-\beta_y\underline{m}'(z)\frac{y_2\left[1+y_2g(z)s(g(z))\right]^3
    h_M(g(z))}
  {1-y_2\int\frac{\left[1+y_2g(z)s(g(z))\right]^2dH(t)}
    {[-t\underline{m}(z)-1-y_2g(z)s(g(z))]^2}}\\
  &=&\frac{\beta_y\underline{m}'(z)y_2\left[1+y_2g(z)s(g(z))\right]^3}
  {1-y_2\int\frac{\left[1+y_2g(z)s(g(z))\right]^2dH(t)}
    {[-t\underline{m}(z)-1-y_2g(z)s(g(z))]^2}}
  \frac{1}{\underline{m}^3(z)}\int\frac{t}{\left(t+m_0(z)\right)^3}dH(t)\\
  &=&-\beta_y\cdot\underline{m}'(z)\frac{y_2\int\frac{tm_0^3(z)}{\left(t+m_0(z)\right)^3}dH(t)}
       {1-y_2\int\frac{m_0^2(z)dH(t)} {(t+m_0(z))^2}}\\
       &=&\frac{\beta_y}{2}\frac{d\left(1-y_2\int\frac{m_0^2(z)dH(t)}
         {(t+m_0(z))^2}\right)}{dz},
\end{eqnarray*}
and
\begin{eqnarray*}
  (\ref{var3})&=&\frac{\beta_yy_2\um'(z_1)\um'(z_2)}{\um^2(z_1)\um^2(z_2)}
  \frac{\partial^2\left[\left(1+y_2g(z_1)
      s(g(z_1))\right)
      \left(1+y_2g(z_2)s(g(z_2))\right)
      h\left(g(z_1),
      g(z_2)\right)\right]}{\partial
    (-1/\underline{m}(z_1))\partial(-1/\underline{m}(z_2))}\\
  &=&\beta_yy_2\frac{\partial^2\int\frac{m_0(z_1)m_0(z_2)}{(t+m_0(z_1))(t+m_0(z_2))}dH(t)}{\partial
    z_1\partial z_2}\\
  &=&\beta_yy_2\int\frac{t^2dH(t)}{(t+m_0(z_1))^2(t+m_0(z_2))^2}\cdot\frac{\partial^2m_0(z_1)m_0(z_2)}{\partial
    z_1\partial z_2}.
\end{eqnarray*}
We have
$$
(\ref{m40})=-\frac{1}{2\pi i}\oint
f(z)\cdot(\ref{m4})dz=\frac{\beta_y}{4\pi i}\oint
f(z)d\left(1-y_2\int\frac{m_0^2(z)dH(t)} {(t+m_0(z))^2}\right).
$$
and
\begin{eqnarray*}
  (\ref{var30})&=&-\frac{1}{4\pi^2}\oint\oint
  f_i(z_1)f_j(z_2)\cdot(\ref{var3})dz_1dz_2\\
  &=&-\frac{\beta_y y_2}{4\pi^2}\oint\oint f_i(z_1)f_j(z_2)
  \left[\int\frac{t^2dH(t)}{(t+m_0(z_1))^2(t+m_0(z_2))^2}\right]dm_0(z_1)dm_0(z_2)~.
\end{eqnarray*}

\subsection{Proof of Lemma \ref{LemCom1}}

Given $z=x_z+i\cdot y_z$, then
\begin{eqnarray*}
  x_z+iy_z&=&-\frac{h^2m_0(z)}{y_2\left(1-y_2+y_2\int\frac{t}{t+m_0(z)}dH(t)\right)}+\frac{y_{1}m_0(z)}{y_{2}}\\
  &=&-\frac{h^2m_0(z)}{y_2\left(1-y_2+y_2\sum\limits_{j=1}^p\frac{ w_j\lambda_j^0 }
    {\lambda_j^0+m_0(z)}\right)}+\frac{y_{1}m_0(z)}{y_{2}}\\
  &=&-\frac{h^2m_0(z)}{y_2\left(1-y_2+y_2\sum\limits_{j=1}^p
    \frac{ w_j \lambda_j^0(\lambda_j^0+u_0 )}{(\lambda_j^0+u_0 )^2+v_0^2 }
    -{\bf i}\cdot y_2\sum\limits_{j=1}^p\frac{ w_j \lambda_j^0 v_0 }{(\lambda_j^0+u_0 )^2+v_0^2 }\right)}+\frac{y_{1}m_0(z)}{y_{2}}\\
  &=&-\frac{h^2 u_0 \left(1-y_2+y_2\sum\limits_{j=1}^p
    \frac{ w_j\lambda_j^0(\lambda_j^0+u_0)}{(\lambda_j^0+u_0)^2+v_0^2 }\right)
    -h^2 v_0 \sum\limits_{j=1}^p\frac{y_2w_j\lambda_j^0v_0}
    {(\lambda_j^0+u_0 )^2+v_0^2}}
  {y_2\left(1-y_2+\sum\limits_{j=1}^p\frac{y_2w_j\lambda_j^0(\lambda_j^0+u_0 )}{(\lambda_j^0+u_0)^2
      +v_0^2 }\right)^2
    +y_2\left(\sum\limits_{j=1}^p\frac{y_2 w_j \lambda_j^0 v_0 }{(1+\lambda_j^0 u_0 )^2
      +(\lambda_j^0)^2 v_0^2 }\right)^2}+\frac{y_1 u_0 }{y_2}\\
  &&{\bf i}\cdot\left[-\frac{h^2 u_0 \sum\limits_{j=1}^p\frac{y_2 w_j \lambda_j^0 v_0 }{(\lambda_j^0+u_0 )^2
        +v_0^2}
      +h^2 v_0 \left(1-y_2+\sum\limits_{j=1}^p\frac{y_2w_j\lambda_j^0(\lambda_j^0+u_0)}{(\lambda_j^0+u_0)^2
        +v_0^2}\right)}
    {y_2\left(1-y_2+\sum\limits_{j=1}^p\frac{y_2 w_j \lambda_j^0(\lambda_j^0+u_0)}{(\lambda_j^0+u_0)^2
        +v_0^2}\right)^2
      +y_2\left(\sum\limits_{j=1}^p\frac{y_2 w_j \lambda_j^0 v_0 }{(\lambda_j^0+u_0 )^2+v_0^2}\right)^2}
    +\frac{y_1 v_0 }{y_2}\right].
\end{eqnarray*}
So we obtain
$$
x_z=-\frac{h^2 u_0 \left(1-y_2+y_2\sum\limits_{j=1}^p
  \frac{ w_j\lambda_j^0(\lambda_j^0+u_0)}{(\lambda_j^0+u_0)^2+v_0^2 }\right)
  -h^2 v_0 \sum\limits_{j=1}^p\frac{y_2w_j\lambda_j^0v_0}
  {(\lambda_j^0+u_0 )^2+v_0^2}}
{y_2\left(1-y_2+\sum\limits_{j=1}^p\frac{y_2w_j\lambda_j^0(\lambda_j^0+u_0 )}{(\lambda_j^0+u_0)^2
    +v_0^2 }\right)^2
  +y_2\left(\sum\limits_{j=1}^p\frac{y_2 w_j \lambda_j^0 v_0 }{(1+\lambda_j^0 u_0 )^2
    +(\lambda_j^0)^2 v_0^2 }\right)^2}+\frac{y_1 u_0 }{y_2}~,
$$
and
$$
y_z=-\frac{h^2 u_0 \sum\limits_{j=1}^p\frac{y_2 w_j \lambda_j^0
    v_0 }{(\lambda_j^0+u_0 )^2
    +v_0^2}
  +h^2 v_0 \left(1-y_2+\sum\limits_{j=1}^p\frac{y_2w_j\lambda_j^0(\lambda_j^0+u_0)}{(\lambda_j^0+u_0)^2
    +v_0^2}\right)}
{y_2\left(1-y_2+\sum\limits_{j=1}^p\frac{y_2 w_j \lambda_j^0(\lambda_j^0+u_0)}{(\lambda_j^0+u_0)^2
    +v_0^2}\right)^2
  +y_2\left(\sum\limits_{j=1}^p\frac{y_2 w_j \lambda_j^0 v_0 }{(\lambda_j^0+u_0 )^2+v_0^2}\right)^2}
+\frac{y_1 v_0 }{y_2}.
$$
So the proof of Lemma \ref{LemCom1} is completed. \eprf


%

\end{document}